\documentclass[12pt]{amsart}

\usepackage{DKstyle}

\newcommand{\vertiii}[1]{{\left\vert\kern-0.25ex\left\vert\kern-0.25ex\left\vert #1
    \right\vert\kern-0.25ex\right\vert\kern-0.25ex\right\vert}}
\theoremstyle{plain}
\makeatletter
\newcommand*{\rom}[1]{\expandafter\@slowromancap\romannumeral #1@}
\makeatother
\newtheorem*{thm*}{Theorem}

\newcommand{\ind}{\mathrm{Ind}}

\usepackage{caption}
\usepackage[labelfont=rm]{subcaption}
\usepackage{bm}
\usepackage[pagebackref=true, colorlinks]{hyperref}

\hypersetup{pdffitwindow=true,linkcolor=blue,citecolor=blue,urlcolor=blue,menucolor=blue}

\usepackage{comment}


\date{\today}%
\dedicatory{}%
\commby{}%

\title{Effective density for inhomogeneous quadratic forms II: fixed forms and generic shifts}
\author{Anish Ghosh}
\address{School of Mathematics, Tata Institute of Fundamental Research, Homi Bhabha Road, Colaba, Mumbai 400005, India}
\email{ghosh@math.tifr.res.in}
\author{Dubi Kelmer}
\address{Department of Mathematics, Boston College, Chestnut Hill MA 02467-3806, USA}
\email{kelmer@bc.edu}
\author{Shucheng Yu}
\address{Department of Mathematics, Technion, Haifa, Israel}
\email{yushucheng@campus.technion.ac.il}
\thanks{AG was supported by a Government of India, Department of Science and Technology, Swarnajayanti fellowship, a CEFIPRA grant, a MATRICS grant and a grant from the Infosys foundation. AG also acknowledges support of the Department of Atomic Energy, Government of India, under project $12-R\&D-TFR-5.01-0500$. DK and SY were partially supported by NSF CAREER grant DMS-1651563. SY acknowledges that this project has received funding from the European Research Council (ERC) under the European Union's Horizon 2020 research and innovation program (grant agreement No.\ 754475).}
\begin{document}
\begin{abstract}
We establish effective versions of Oppenheim's conjecture for generic inhomogeneous quadratic forms. We prove such results for fixed quadratic forms and generic shifts. Our results complement our previous paper \cite{GhoshKelmerYu2019a} where we considered generic forms and fixed shifts. In this paper, we use ergodic theorems and in particular we establish a strong spectral gap with effective bounds for some representations of orthogonal groups which do not possess Kazhdan's property $(T)$.  
\end{abstract}
\maketitle

\tableofcontents
\section{Introduction}
Let $Q$ be a quadratic form on $\R^{n}$ and let ${\alpha}$ be a vector in $\R^{n}$. Define the inhomogeneous quadratic form $Q_{{\alpha}}$ by
$$Q_{{\alpha}}({v})=Q({v}+{\alpha})\ \textrm{for any}\ {v}\in \R^{n},$$
where we think of $Q_{\alpha}$ as a \textit{shift} by ${ \alpha}$ of the homogenous form $Q$.
The inhomogeneous form $Q_{{\alpha}}$ is said to be \textit{indefinite} if $Q$ is indefinite and \textit{non-degenerate} if $Q$ is non-degenerate. Finally, $Q_{{\alpha}}$ is said to be \textit{irrational} if either $Q$ is an irrational quadratic form, i.e. not proportional to a quadratic form with integer coefficients, or ${\alpha}$ is an irrational vector. 

The famous Oppenheim conjecture admits a natural variant for inhomogeneous forms. Namely, it follows from the work of Margulis and Mohammadi  
\cite{MargulisMohammadi11} (who obtained a more quantitative result on the density of such values) that for any indefinite, irrational,  non-degenerate inhomogeneous form $Q_{ \alpha}$ in $n \geq 3$ variables, $Q_{{\alpha}}(\mathbb{Z}^{n})$ is dense in $\mathbb{R}$. However, inhomogeneous quadratic forms have attracted considerable attention earlier; we refer the reader to J. Marklof's important works \cite{Marklof2002, Marklof2003} on their pair correlation density for example. 
In this paper we are concerned with the question of \emph{effectivity}, namely, for a given $Q_{{\alpha}}, \xi \in \R$ and $t \geq 1$ large: How small can $|Q_{{\alpha}}({v}) - \xi|$ get for  ${v} \in \Z^{n}$ with $\|{v}\|\leq t$ bounded?  This is a notoriously difficult problem intimately connected with questions of effectivity in homogeneous dynamics and has received considerable attention recently. We refer the reader to \cite{GhoshKelmerYu2019a} for a relatively comprehensive account of the work that has been done on homogeneous forms \cite{GotzeMargulis2010, LindenstraussMargulis2014, Bourgain2016, GhoshGorodnikNevo2020, AthreyaMargulis2018, GhoshKelmer17, GhoshKelmer2018, KelmerYu2020, KleinbockSkenderi2021} and on inhomogeneous forms \cite{StrombergssonVishe2020}. 

In a previous paper \cite{GhoshKelmerYu2019a} we considered this question for generic inhomogeneous forms. There is a natural measure on the space of forms (of a fixed signature and discriminant) and using a second moment formula for Siegel transforms, it can be shown \cite[Theorem 1.1]{GhoshKelmerYu2019a}, that for any $\kappa<n-2,$ almost all indefinite forms in $n$ variables, and almost all shifts $\alpha\in \R^n,$ the system of inequalities 
\begin{equation} \label{equ:maininequ} 
|Q_{{\alpha}}({v})-\xi|<t^{-\kappa},\; \|{v}\|\leq t
\end{equation}
has integer solutions for all sufficiently large $t$. The main result of \cite{GhoshKelmerYu2019a}  addresses the much more difficult problem of effectivity for fixed shifts and generic forms. More precisely,  by proving a second moment formula for congruence groups, we showed that the same result as above holds for any fixed rational $\alpha\in \Q^n$ and almost all indefinite forms. In fact, we obtain a counting result \cite[Theorem 1.2]{GhoshKelmerYu2019a} when the shift is rational. For fixed irrational shifts, we obtain weaker bounds.

In this paper we study the complementary problem of a fixed indefinite form $Q$, and almost all shifts $\alpha$. For this problem we have the following result for rational forms.
\begin{Thm}\label{thm:genericshifts}
For any rational indefinite form $Q$ in $n$ variables and any $\xi\in \R,$
there is $\kappa_0$ $($depending only on the signature of $Q$$)$ such that for any $\kappa<\kappa_0,$  for almost all $\alpha\in \R^n$ the system of inequalities \eqref{equ:maininequ} 
has integer solutions for all sufficiently large $t$.
\end{Thm}
Our proof gives the following explicit values for $\kappa_0$, depending only on the signature of  $Q$ 
$$\kappa_0= \left\lbrace\begin{array}{ll} 
1 & (p,q)=(2,1)\\
2 & (p,q)=(n-1,1),\;n\geq 3\\
2& (p,q)=(2,2)\\
3/2 & (p,q)=(4,2) \mbox{ or } (3,3)\\
5/2 & (p,q)=(6,3),
\end{array}\right.
$$
while for all other signatures $p\geq q>1$ with $p+q=n$ we have $\kappa_0=2\kappa_1q(p-1)$ with 
$$\kappa_1= \left\lbrace\begin{array}{ll} 
\tfrac{1}{n} & n \equiv0\pmod{4}\\
\frac{1}{n-1} & n\equiv 1\pmod{4}\\
\tfrac{1}{n-2}& n\equiv 2\pmod{4}\\
\tfrac{1}{n+1} & n\equiv 3\pmod{4}.
\end{array}\right.
$$
\begin{rem}
For $n=3,4$, our method gives the optimal bound $\kappa_0=n-2$. However, in general when $n\geq 5$ our result is not as good as what we have for generic forms. For example for signature $(n-1,1)$ we have $\kappa_0=2$ which is much smaller than the expected rate of $n-2$. For large values of $n$ our bounds are best when $p= q$ (or $p=q+1$) in which case we get $\kappa_0$ is roughly $n/2$ which is still about half the expected value.
\end{rem}
The method we use for this problem is completely different from the one used to study generic forms. Here we reduce the problem to a \emph{shrinking target problem} for the action of a semi-simple group acting on a homogeneous space, and then rely on an effective mean ergodic theorem to study the shrinking target problem. This is similar to the approach taken in \cite{GhoshGorodnikNevo2020} and also in \cite{GhoshKelmer17,GhoshKelmer2018}.
In order to outline the general idea and also explain where the exponents are coming from we need to introduce some more notation.

For $Q$ an indefinite quadratic form, let $G=\SO^+_Q(\R)$ denote the connected component of the identity in the group of linear transformations preserving $Q$, and note that for a rational form $Q$ we have that the set of integer points $\G=\SO^+_Q(\Z)$ is a lattice in $G$. Using the natural embedding of $G$ in $\SL_n(\R)$ we get a natural action of $G$ on $\R^n$ and we may consider the semi-direct product 
$\tilde{G}=G\ltimes \R^n$.  We note that  $\tilde{\G}=\G\ltimes \Z^n$ is a lattice in $\tilde{G}$ and that there is a natural left action of $G$ on the space $L^2(\tilde{G}/\tilde{\G})$ preserving the probability Haar measure $m_{\tilde{G}}$ on $\tilde{G}/\tilde{\G}$. Our first step is to reduce the problem of approximating a target $\xi$ by values $Q( v+\alpha)$ to a shrinking target problem for the action of $G$ on $\tilde{G}/\tilde{\G}$. 
Using the results of \cite{GhoshKelmer17} we can further reduce this problem to obtaining an appropriate effective mean ergodic theorem. 

For any $f\in L^2(\tilde{G}/\tilde{\G})$ and growing measurable subsets $G_t\subseteq G$ consider the averaging operator
\begin{equation}\label{equ:avgop}
\beta_{G_t}f(x)=\frac{1}{m_G(G_t)}\int_{G_t}f(g^{-1}x)dm_G(g),
\end{equation}
where $m_G$ denotes the probability measure on $G/\G$ coming from Haar measure of $G$. From \cite{GhoshKelmer17}, having an effective mean ergodic theorem of the form 
$$\left\|\beta_{G_t}f-\int_{\tilde{G}/\tilde{\G}}fdm_{\tilde{G}} \right\|_2 \leq C \frac{\|f\|_2}{m_G(G_t)^\kappa},$$
valid for all $f\in L^2(\tilde{G}/\tilde{\G})$ has very strong consequences for shrinking target problems for the action of $G$ on $\tilde{G}/\tilde{\G}$. 
In fact, we will show that in our case, it is enough to have such a result for functions in the smaller space $L^2_{00}(\tilde{G}/\tilde{\G})$ of all functions whose average over $\R^n/\Z^n$ is zero. As will become apparent below, this refinement is crucial in order to get uniform bounds for signatures $(n-1,1)$ and $(2,2)$. More precisely, we will show that  \thmref{thm:genericshifts} follows from the following.

%We can naturally identify $L^2(G/\G)$ as a subspace of $L^2(\tilde{G}/\tilde{\G})$ and denote by $L^2_0(\tilde{G}/\tilde{\G})$ its orthogonal complement. We thus have a decomposition 
%$$L^2(\tilde{G}/\tilde{\G})=L^2(G/\G)\oplus L^2_0(\tilde{G}/\tilde{\G}),$$
%which is preserved by the $G$ action. Let $G_t$ be a family of growing norm balls in $G$ given by 
%$$G_t=\{g\in G\ |\ \|g\|_{\mathrm{op}}\leq t\},$$
%where $\|\cdot\|_{\mathrm{op}}$ is the operator norm on $\SL_n(\R)$ with respect to the Euclidean norm on $\R^n$. Let $\pi$ denote the representation of $G$ on the space $L^2_0(\tilde{G}/\tilde{\G})$ and consider the averaging operators
%$$\beta_{G_t}(\pi)(v)=\int_{G_t}\pi(g)vdm(g).$$
%We will then show the following effective mean ergodic theorem 

\begin{Thm}\label{t:EMET}
Let $Q$ be an indefinite rational form of signature $(n-q,q)$ and let $\kappa_1=\frac{\kappa_0}{2q(n-q-1)}$ with $\kappa_0$ as above. There is a family of growing norm balls $G_t\subseteq \SO^+_Q(\R)$ of measure $m_G(G_t)\gg t^{q(n-q-1)}$ such that for any $\kappa<\kappa_1$ and for any $f\in L^2_{00}(\tilde{G}/\tilde{\G})$ % with $\int_{\R^n/\Z^n}f(g,v)dv=0,$ 
we have that  
$$\|\beta_{G_t} f\|_2\ll_{\kappa} \frac{\|f\|_2}{m_G(G_t)^\kappa},$$
where the implied constant depends only on $\kappa$.
\end{Thm}
%We will show that \thmref{thm:genericshifts} for a rational form of signature $(p,q)$ with $p\geq q\geq 1$ follows from \thmref{t:EMET} with the critical exponent given by $\kappa_0=2\kappa_1q(p-1)$. 

For the proof of \thmref{t:EMET} we will exploit a general spectral transfer principle, described in \cite{Nevo1998}, giving  explicit bounds on the exponent $\kappa_1$ in the mean ergodic theorem in terms of the strong spectral gap of the corresponding representation. This  further reduces our problem, to obtaining effective bounds for the strong spectral gap for the representation of $G$ on the space $L^2_{00}(\tilde{G}/\tilde{\Gamma})$, which is a problem of independent interest.

The strong spectral gap of a unitary representation $\pi$ is controlled by the parameter $p(\pi)\geq 2$ which can be defined as follows.
We say that a unitary representation $\pi$ of $G$ on a Hilbert space $\cH$ is \textit{strongly $L^p$} if the functions $g\mapsto \langle \pi(p)v,v\rangle$ are all in $L^p(G)$ for a dense set of vectors $v$ in $\cH$. The parameter $p(\pi)$ is then the infimum over all $p\geq 2$ for which $\pi$ is strongly $L^p$. We note that a unitary representation $\pi$ is called \textit{tempered} if $p(\pi)=2$.
% The problem of obtaining bounds for the strong spectral gap is a problem of independent 
% We recall that the strong spectral gap of a unitary representation $\pi$ of a semi-simple connected group $G$, measures how far it is from being tempered, and can be captured by the parameter, $p(\pi)$ we now define.  defined as the infimum over all $p\geq 2$ for which the functions $g\mapsto \langle \pi(p)v,u\rangle$ are all in $L^p(G)$ for a dense collection of vectors in $\pi$. Using this notion one can get a bound for the critical exponent $\kappa_1$  of the order of $1/p(\pi)$ (see  \secref{s:OperatorNorm} or precise statement). 

%For $n=3$, it follows from a classical result of Kazhdan \cite{Kazhdan1967} that any representation of $\SO^+_Q(\R)\ltimes \R^3$ with no nontrivial $\R^3$-invariant vectors is tempered. Hence for $n=3$ we get the optimal exponent $\kappa_0=1$. In general we do not know if these representations are necessarily tempered (and as we show below they are in general non-tempered).
%%\comm{(comm: note that whenever $\kappa>n-2$ then it follows from \cite{KleinbockSkenderi} that \eqref{equ:maininequ} with $\xi=0$ has no integer solutions for all sufficiently large $t$, so we can actually conclude here that if $q(p-1)>n-2$ then the representation on $L^2_{00}(\tilde{G}/\tilde{\G})$ is not tempered)}. 
%Nevertheless, 

When the group $G$ has Kazhdan's property $(T)$ there is a uniform bound $p(G)$ such that $p(\pi)\leq p(G)$ for all nontrivial irreducible unitary representations of $G$. Moreover, it was shown by Wang \cite{Wang2014}, that %for any simple linear group $G$ acting linearly on $\R^n$ with no nontrivial fixed points (denote this action by $\rho$), 
there is a constant $p(\tilde{G})_0$ such that $p(\pi)\leq p(\tilde{G})_0$ for any representation $\pi=\tilde{\pi}|_G$, with $\tilde{\pi}$ a unitary representation of $\tilde{G}=G\ltimes \R^n$, with no nontrivial $\R^n$-invariant vectors. (Indeed Wang studied this question in a much greater generality of a simple linear group acting linearly on a Euclidean space with no nontrivial fixed points.)
%Specifying to our case when $G=\SO_Q^+(\R)$ (acting on $\R^n$ via the natural embedding of $G$ in $\SL_n(\R)$) 

For any finite-index subgroup $\Gamma\leq \SO_Q^+(\Z)$ we let $\pi_{\tilde\Gamma}$ denote the representation of $G$ on $L^2_{00}(\tilde{G}/\tilde\Gamma)$ and let, $p(\tilde{\G})_0=p(\pi_{\tilde\Gamma})$. Since $\pi_{\tilde{\Gamma}}$ is the restriction of a representation of $\tilde{G}$ with no nontrivial $\R^n$-invariant vectors, we get the hierarchy 
$$p(\tilde{\Gamma})_0\leq p(\tilde{G})_0\leq p(G).$$

When $G$ has property $(T)$ the work of \cite{Li1995,Oh98} give very good bounds for $p(G)$ (which are shown in many cases to be sharp). However, much less is known regarding the values of $p(\tilde{G})_0$ and $p(\tilde\Gamma)_0$. For the special case when $G= \SO_Q^+(\R)$ with $Q$ of signature $(2,1)$, it follows from a classical result of Kazhdan \cite{Kazhdan1967} that any unitary representation of $\tilde{G}=\SO^+_Q(\R)\ltimes \R^3$ with no nontrivial $\R^3$-invariant vectors is tempered, implying the optimal bounds  $p(\tilde{G})_0=2$.   Going beyond this case, the method of Wang \cite{Wang2014}, produces explicit upper bounds for $p(\tilde{G})_0$ in quite great generality. However, these bounds are usually far from optimal (in particular, when $G$ has property $(T)$ the known bounds for $p(G)$ are usually better). It is thus an interesting problem to give sharp bounds for $p(\tilde{G})_0$ and $p(\tilde\Gamma)_0$.

Our first result in this regard is the following lower bound for $p(\tilde{\G})_0$ (and hence also for $p(\tilde{G})_0$), for $G=\SO_Q^+(\R)$ with $Q$ an indefinite rational form and $\G=\SO^+_Q(\Z)$,  
showing that in most cases these representations are not tempered.
 \begin{Thm}\label{t:LowerBound}
For $Q$ an indefinite rational form of signature $(n-q,q)$ with $2q\leq n$, let $G=\SO_Q^+(\R)$ and $\Gamma=\SO^+_Q(\Z)$ and let $\tilde{G}=G\ltimes \R^n$ and $\tilde\G=\G\ltimes \Z^n$ be as above. With the exception of signatures $(2,1),(3,1), (2,2)$ and possibly $(3,2)$ the representation of $G$ on  $L^2_{00}(\tilde{G}/\tilde{\G})$ is non-tempered.  Moreover, for $Q$ of signature $(n-1,1)$ with $n\geq 4$ we have $p(\tilde{\Gamma})_0\geq n-2$, for $Q$ of signature $(n-q,q)$ with $1<q<\tfrac{n}{2}$ we have 
$p(\tilde{\G})_0> n-q-2-j$ with $j\in\{0,1,2,3\}$ satisfying $j\equiv n-q-2\pmod{4}$ and for signature $(q,q)$ we have $p(\tilde{\G})_0> n-q-1-j$  with $j\in\{0,1,2,3\}$ satisfying $j\equiv n-q-1\pmod{4}$.
\end{Thm}

In terms of upper bounds, when $G$ has property $(T)$ we could not improve over the general bounds for $p(G)$. 
Our new results in this regard are for forms of signature $(2,2)$ and $(n-1,1)$. Our first result is an optimal upper bound for forms with signature $(2,2)$, analogous to Kazhdan's result for signature $(2,1)$.

\begin{Thm}\label{t:SG}
Let $Q$ be a rational form of signature $(2,2)$, let $G=\SO_Q^+(\R)$ and let $\tilde{G}=G\ltimes \R^4$. Let $\pi=\tilde{\pi}|_{G}$ with $\tilde{\pi}$ a unitary representation of $\tilde{G}$ without nontrivial  $\R^4$-invariant vectors, then $\pi$ is tempered.
\end{Thm}

Next we consider a form $Q$ of signature $(n-1,1)$ and as before  $G=\SO_Q^+(\R)$ and $\tilde{G}=G \ltimes \R^n$. Here, a refinement of the general argument of Wang (see \cite[Example 4 ]{Wang2014}) implies that $p(\tilde{G})_0\leq 2(n-2)$. We conjecture that the correct bound in this case is $p(\tilde{G})_0= (n-2)$. As evidence for this conjecture we prove it for the representations of $G$ on $L^2_{00}(\tilde{G}/\tilde{\G})$ (which is what we need for our application). 
 
\begin{Thm}\label{t:rank1a}
Let $Q$ be a rational form of signature $(n-1,1)$  with $n\geq 4$, let $G=\SO^+_Q(\R)$ and let $\tilde{G}=G\ltimes \R^n$ be as above. For any finite-index subgroup $\G\leq \SO^+_Q(\Z)$ 
we have $p(\tilde\Gamma)_0=n-2$. 
\end{Thm}
%\comm{[Looking at the proof of this theorem, we chose specially the standard $(n-1,1)$ form for some of our computation; is it clear we can translate this result to all other $(n-1,1)$ forms? It seems for $G=h^{-1}\SO_{Q_0}^+h$ we can choose our Cartan decomposition $G=KA^+K$ with $K=h^{-1}K_0h$ and $A=h^{-1}A_0k$, then all the relations still stay true.]}
\begin{rem}
It is interesting to compare this result to the analogous  result on the spectral gap $p(\Gamma)$ of the representation of $G$ on $L^2(G/\Gamma)$.
When $\Gamma$ is a congruence group, the  Selberg-Ramanujan conjecture states that $p(\Gamma)\leq \max\{2, n-2\}$ for any congruence group $\Gamma\leq \SO^+_Q(\Z)$, whereas the best known bounds towards this conjecture are currently $p(\Gamma)\leq \frac{64}{25}$ for $n=3,4$ \cite{KimSarnak03,BlomerBrumley11}, while for $n\geq 5$ the conjectured bound was proved in \cite{BergeronClozel13}. These bounds no longer hold for general finite-index subgroups $\Gamma\leq \SO^+_Q(\Z)$, as there are constructions of finite-index subgroups $\Gamma\leq \SO^+_Q(\Z)$ with $p(\Gamma)$ arbitrarily large. This is very different from our current case, where the bound for $p(\tilde\Gamma)_0$ holds for all finite-index subgroups (which is compatible with our conjecture that $p(\tilde{G})_0=n-2$ in this case).
\end{rem}
Finally, as further evidence for our conjecture we prove the following theorem reducing the conjecture to the case of signature $(3,1)$.
\begin{Thm}\label{t:rank1b}
If $p(\tilde{G})_0=2$ for $G\cong\SO^+(3,1)$ then $p(\tilde{G})_{0}=n-2$ for $G\cong\SO^+(n-1,1)$ for all $n\geq 5$.
\end{Thm}

\subsection*{Acknowledgements}
We thank Alireza Salehi Golsefidy, Amir Mohammadi and Amos Nevo for our discussions on this problem.

\section{Preliminaries and notation}
\subsection{Notation}
Let $n=p+q\geq 3$ with $p\geq q\geq 1$. Let $Q$ be a non-degenerate quadratic form of signature $(p,q)$. Then $Q$ can be represented by a unique invertible symmetric matrix $J\in M_n(\R)$ in the sense that $Q(v)=vJv^*$ for any $v\in \R^n$, where $v^*$ denotes the transpose of $v$. Let $G=\SO^+_Q(\R)$ denote the connected component of the identity inside the special orthogonal group preserving $Q$. We use the notation $A\ll B$ as well as $A=O(B)$ to indicate that there is some constant $c>0$ such that $A\leq cB$. The constant may depend on $n$ that we think of as fixed, if we want to emphasize the dependance of the constant on various parameters we will indicate it with a subscript. We also use the notation $A\asymp B$ to mean that $A\ll B\ll A$.
\subsection{Coordinates}\label{subsec:coordinate}
For some calculations we will need to work with explicit coordinates. Since, for any form of signature $(p,q)$ the group 
$\SO_Q(\R)$ can be conjugated in $\SL_n(\R)$ to the group $\SO_{Q_0}(\R)$ with 
\begin{equation}\label{equ:baseform}
Q_0(v)=\sum_{i=1}^pv_i^2-\sum_{i=p+1}^nv_i^2,
\end{equation}  
it is enough to consider the case of $Q=Q_0$.
The group $G$ has a polar decomposition $G=KA^+K$ with $K$ a maximal compact subgroup and $A^+$ the positive Weyl chamber in the Cartan group $A$. Explicitly, for $G=\SO_{Q_0}^+(\R)$ we can take
the maximal compact subgroup
$$K=\left\{k=\begin{pmatrix}
k_1& \\
 & k_2\end{pmatrix}\ |\ k_1\in \SO_p(\R)\ \textrm{and}\ k_2\in \SO_q(\R)\right\}.$$
and the Cartan group $A=\exp{\fa}$ with 
$$\fa=\{H=\diag^{\star}\ (h_1,\ldots, h_q,0,\ldots, 0, h_q,\ldots,h_1)\ |\ h_i\in \R\ \textrm{for}\ 1\leq i\leq q\},$$
where $\diag^{\star}(t_1, t_2, \ldots,t_n)$ denotes the anti-diagonal matrix with the $(n-i,i)^{\textrm{th}}$ entry given by $t_i$. 
Moreover, the positive Weyl chamber $\fa^+$ can be taken such that if $p>q$ then $H\in \fa^+$ if $h_1\geq h_2\geq\cdots\geq h_q\geq 0$, and if $p=q$ then $H\in \fa^+$ if $h_1\geq h_2\geq \cdots\geq h_{q-1}\geq |h_q|$.

We can now describe the Haar measure $m_G$ of $G$ as follows (see \cite[Proposition 5.28]{Knapp1986}): For any $g\in G$ writing $g=k_1\exp(H)k_2$ with $k_1,k_2\in K$ and $H\in \fa^+$, then up to a normalizing factor
\begin{equation}\label{equ:Haar}
dm_G(g)=\prod_{1\leq i<j\leq q}\sinh(h_i-h_j)\sinh(h_i+h_j)\prod_{1\leq i\leq q}\sinh(h_i)^{p-q}dk_1dHdk_2,
\end{equation}
where $dk$ is the probability Haar measure on $K\cong \SO_p(\R)\times \SO_q(\R)$ and $dH$ is the Lebesgue measure on $\fa$ (identified with $\R^q$).

For the case of signature $(n-1,1)$ we will need to make some explicit calculations, so we give some more details on this decomposition. In this case we have the polar decomposition $G=KA^+K$ with
$$K=\left\{k=\left(\begin{smallmatrix}
k' & \\
 & 1\end{smallmatrix}\right)\ |\ k'\in \SO_{n-1}(\R)\right\}\ \textrm{and}\ A^+=\left\{a_t=\left(\begin{smallmatrix}
 \cosh t & & \sinh t\\
  &I_{n-2}& \\
  \sinh t & &\cosh t\end{smallmatrix}\right)\ |\ t>0\right\}.$$
Let $M\cong \SO_{n-2}(\R)$ be the centralizer of $A$ in $K$, namely
$$M=\left\{m=\left(\begin{smallmatrix}
1 & & \\
 &m'& \\
 &   &1\end{smallmatrix}\right)\ |\ m'\in \SO_{n-2}(\R)\right\}.$$
 It was shown in \cite[Lemma 2.1 and Proposition 2.4]{RisagerSodergen16} that any $k\in K$ can be written as 
 $k=m_1k_{\theta}m_2$, where $m_i\in M$ and 
\begin{equation}\label{e:ktheta}
k_{\theta}=\left(\begin{smallmatrix}
\cos\theta & \sin\theta & \\
-\sin\theta & \cos\theta & \\
 & & I_{n-2}\end{smallmatrix}\right)
 \end{equation}
for some uniquely determined $\theta\in [0,\pi]$, and that with this decomposition, the Haar measure $dk$ of $K$ is given (up to a normalizing constant) by
\begin{equation}\label{e:HaarMKM}
dk=\left(\sin\theta\right)^{n-3} dm_1d\theta dm_2.
\end{equation}

%\begin{remark}
%For a general rational form of signature $Q$ of signature $(p,q)$, there is $h\in \SL_n(\R)$ such that $G'=\SO_Q^+(\R)=h^{-1}\SO_{Q_0}^+(\R)h$ 
%\end{remark}

\subsection{Norm balls}
We denote by $\|\cdot\|$ the Euclidean norm on $\R^n$ and using the natural embedding of $G=\SO^+_Q(\R)\subseteq \SL_n(\R)\subseteq \GL_n(\R)$ we let $\|\cdot\|_{\rm op}$ denote the operator norm on $G$. Using the operator norm we define the norm $\|g\|=\|g^{-1}\|_{\rm op}$ and we use it to define the growing norm balls  
\begin{equation}\label{e:G_T}
G_T=\{g\in G\ |\  \|g\|\leq T\}.
\end{equation}
\begin{rem}\label{rmk:comofnorm}
This choice of norm is convenient for what follows but is not essential. Note that for $g\in \SO_Q(\R)$ we have that  $g^{-1}=Jg^*J^{-1}$ implying that $\|g^{-1}\|_{\rm op}\asymp \|g\|_{\rm op}$. 
%For later use, we fix $\tilde{C}>1$ depending only on $Q$ such that 
%\begin{align}\label{equ:comnorm}
%\tilde{C}^{-1}\|g^{-1}\|\leq \|g\|\leq \tilde{C}\|g^{-1}\|,\quad \forall\ g\in G.
%\end{align}
We can also use the Hilbert-Schmidt norm instead of the operator norm, noting that both the Hilbert-Schmidt norm and the operator norm are bi-$\operatorname{O}_n(\R)$-invariant and that $\|g\|_{\rm HS}\asymp \|g\|_{\rm op}$ for all $g\in \GL_n(\R)$. 
\end{rem}

The volume of such norm balls was computed (in greater generality) in \cite[
Corollary 1]{GorodnikWeiss2007} and satisfy that  there is $T_0$ depending only on $n$ such that for all $T\geq T_0$
\begin{equation}\label{e:mG_T}
m_G(G_T)\asymp (\log T)^{\delta_{p,q}}T^{q(p-1)}.
\end{equation}
We also note that since the identity matrix satisfies $ \| I \| =1$, from the continuity of the norm for any $T>1,$ the set $G_T$ contains some neighborhood of the identity, and hence $m_G(G_T)>0$ for all $T>1$. Finally for later use we note that when $G=\SO_{Q_0}^+(\R)$ with $Q_0$ of signature $(n-1,1)$, we have that for any $T>1$
\begin{equation}\label{equ:gtp1}
G_T=\{k_1a_tk_2\ |\ k_1,k_2\in K\ \textrm{and}\ 0<t\leq \log T\}.
\end{equation}

\section{Reduction to an effective mean ergodic theorem}\label{sec:reduction}
In this section we will perform several reductions to the problem and we will assume throughout this section that $Q$ is a rational form so that $\G=\SO^+_{Q}(\Z)$ is a lattice in $G=\SO^+_Q(\R)$. First we reduce the problem to a shrinking target problem for the action of $\G=\SO^+_Q(\Z)$ on the torus $\bT^n=\R^n/\Z^n$ (such problems were studied in detail in \cite{GhoshGorodnikNevo14}). Next, we further reduce the problem to a shrinking target problem for the action of $G$ on the space $G \ltimes \R^n/\G\ltimes \Z^n$. Finally we show how this second shrinking target problem follows from an appropriate effective mean ergodic theorem.

\subsection{Reduction to a shrinking target problem for $\G$-action}
Since $\G=\SO^+_Q(\Z)$ is naturally embedded in $\SL_n(\Z),$ it acts on the torus $\bT^n=\R^n/\Z^n$. We call a family, $\{A_t\}_{t>0}$, of subsets of $\bT^n$ a family of \textit{shrinking targets} if $A_t\subseteq A_s$ for all $t\geq s$, $m(A_t)>0$ for all $t>0$ and $m(A_t)\to 0$ as $t\to \infty$ (where $dm$ stands for Lebesgue's measure). 
\begin{Lem}\label{lem:Reduction1}
Assume that there exists some $\kappa_0>0$ such that for any $\kappa<\kappa_0,$ for any family of shrinking targets $\{A_t\}_{t>0}$ in $\bT^n$ with $m(A_t)\gg t^{-\kappa},$ for almost every $\alpha\in \bT^n$ and for all sufficiently large $t,$ there is $\g\in\G$ with $\|\g\|\leq ct$ and $\alpha \gamma\in A_t$ for some $c>0$. Then for any $\kappa<\kappa_0$,  for almost every $\alpha\in\R^n$ and for all sufficiently large $t,$ there is $ v\in \Z^n$ with $|Q( v+\alpha)-\xi|< t^{-\kappa}$ and $\|v\|\leq t$.
\end{Lem} 
\begin{proof}
Fix $\kappa<\kappa_0$ and let $\kappa'\in(\kappa,\kappa_0)$.  
Let
$$\tilde{A}_{N,\e}=\{x \in \R^n\ |\ \|x\|\leq N,\; |Q(x)-\xi|\leq \e\}.$$
Then by \cite[Theorem 5]{KelmerYu2020} we have
$m(A_{N,\e})=2c_Q\e N^{n-2}(1+O_Q(N^{-1/2}))$ and in particular there is some $N_0$  (depending only on $Q$ and $\xi$) such that $m( A_{N_0,\e}) \asymp_{N_0} \e$ for all $\e\in(0,1)$.
Since $N_0$ is fixed, the set $\tilde A_{N_0,\e}$ is contained in finitely many fundamental domains for $\bT^n$ (with the number of fundamental domains uniform in $\e\in (0,1)$); so denoting by $A_t$ the projection of $\tilde{A}_{N_0,t^{-\kappa'}}$ to $\bT^n,$ we still have that $m(A_t)\asymp t^{-\kappa'}$.

Now by assumption, for almost every $\alpha\in \bT^n$ and for all sufficiently large $t$, there is $\g\in \G$ with 
$\|\g\|\leq ct$ and $\alpha\g\in A_t$. Hence for almost every $\alpha\in \R^n$ and for all sufficiently large $t$, there is $\g\in \G$ with $\|\g\|\leq ct$ and $u\in \Z^n$ with $ x=\alpha\gamma+ u\in \tilde{A}_{N_0,t^{-\kappa'}}$. Let $ v= u\gamma^{-1}\in \Z^n$ then since $\g\in \G=\SO^+_Q(\Z),$ we have that
$$Q( v+\alpha)=Q( v\gamma+\alpha\gamma)=Q(x),$$
so that $|Q( v+ \alpha)-\xi|<t^{-\kappa'}$. 
%Finally by \lemref{lem:inversenorm} there is an absolute constant $C>0$ satisfying that $\|g^{-1}\|\leq C\|g\|$ for any $g\in G$. 
We can thus estimate for all $t$ sufficiently large
$$\| v\|=\| u\gamma^{-1}\|=\| x\g^{-1}+\alpha\|\leq N_0\|\gamma^{-1}\|_{ \rm op}+\|\alpha\|\leq N_0ct+1\leq 2cN_0t.$$
Hence, replacing $2cN_0t$ by $t$ we get that for almost every $\alpha\in \R^n$ for all sufficiently large $t$ there is $ v\in \Z^n$ with $\| v\|\leq t$ 
and $|Q( v+ \alpha)-\xi|\leq (t/2cN_0)^{-\kappa'}<t^{-\kappa}$.
\end{proof}

\subsection{Reduction to shrinking target problem for $G$-action}
We now show that the shrinking target property we need for the $\G$-action follows from an appropriate shrinking target property for the action of $G=\SO^+_Q(\R)$ on the space $\tilde{G}/\tilde{\G}$ 
with $\tilde{G}=G\ltimes \R^n$ and $\tilde{\G}=\G\ltimes \Z^n$ (here the group law on $\tilde{G}$ is given by $(h,\alpha)(g,v)=(hg,\alpha g+v)$ for any $(h,\alpha),(g,v)\in \tilde{G}$). %Explicitly, we show the following.

\begin{Lem}\label{lem:Reduction2}
Let $A_t\subseteq \bT^n$ be a family of shrinking targets and let $\tilde{A_t}\subseteq [0,1]^n$ be sets whose projection to $\bT^n$ equals $A_t$. Let $\cF\subset G$ be a fixed fundamental domain for $G/\G$ containing the identity element. Fix a small constant $c>1$ and consider the sets  $\tilde{B_t}=\{(g,\alpha)\in \tilde{G}\ |\ g\in \cF, \|g\|\leq c, \alpha\in \tilde{A}_t\}$ in $\tilde{G}$, 
%where $c>1$ is sufficiently large so that $\tilde{B_t}$ is of positive measure with respective to some Haar measure of $\tilde{G}$, 
and let $B_t\subseteq \tilde{G}/\tilde{\G}$ denote the projection of $\tilde{B_t}$ to $\tilde{G}/\tilde{\G}$. If for a.e. $x\in \tilde{G}/\tilde{\G}$ and for all sufficiently large $t$, there is $g\in G_t$ with $g^{-1}x\in B_t$ then for a.e. $\alpha\in \bT^n$ and all sufficiently large $t$, there is $\g\in \G$ with $\|\g\|\leq c^2t$ and $\alpha \g\in A_t$.
\end{Lem}
\begin{proof}
Since we assume that the set of $x\in \tilde{G}/\tilde{\G}$ such that for all sufficiently large $t$ there is $g\in G_t$ with $g^{-1}x\in B_t$ is a set of full measure, by unfolding we also have that the set
$$S=\{(h,\alpha)\in \tilde G\ |\ \exists\ t_0> 1\ s.t.\ \forall\ t\geq t_0\  \exists\ g\in G_t, \g\in \G, m\in \Z^n\ s.t.\ (g^{-1},0)(h,\alpha)(\g,m)\in \tilde{B_t}\},$$
is a set of full measure in $\tilde{G}$. In particular, if $\tilde{G}_c =\{(h,\alpha)\in \tilde{G}\ |\ \|h^{-1}\|\leq c\}$ then $S'=S\cap \tilde{G}_c$ is of full measure in $\tilde{G}_c$. Moreover, since $c>1$ the set $\{h\in G\ |\ \|h^{-1}\|\leq c\}$ has positive measure. Hence for almost every $\alpha\in \R^n$ there exists some $h\in G$ with $\|h^{-1}\|\leq c$, for which $(h, \alpha)\in S'$, that is, for all sufficiently large $t$ there is $g\in G_t,\g\in \G$ and $m\in \Z^n$ such that 
$$(g^{-1}h\g,\alpha \g+m)\in \tilde{B_t}.$$ 
For such a pair $(h,\alpha)$ we have that $\|g^{-1}h\g\|\leq c$ and $\|h^{-1}\|\leq c$ so 
$$\|\g\|=\|h^{-1}gg^{-1}h\g\|\leq \|h^{-1}\|\|g\|\|g^{-1}h\g\|\leq c^2t,$$
and that $\alpha \g+m\in \tilde{A_t}$ so that $\alpha \g+\Z^n\in A_t$.
We thus showed that for a.e. $\alpha \in \bT^n$ and for all sufficiently large $t$ there is $\g\in \G$ with $\|\g\|\leq c^2t$ such that $\alpha\g\in A_t$ concluding the proof.
\end{proof}

\subsection{Reduction to an effective mean ergodic theorem}
Given a measure preserving ergodic action of a noncompact, locally compact group $G$ on a probability space $(X,m_X)$, the mean ergodic theorem states that 
for any growing family of subsets $G_t$ of $G$, the averaging operator on $L^2(X)$ given by
\begin{equation}\label{e:AveOp}
\beta_{G_t}f(x)=\frac{1}{m_G(G_t)}\int_{G_t}f(g^{-1} x)dm_G(g),
\end{equation}
satisfies that $\|\beta_{G_t}f-\int_{X}fdm_X\|_2\to 0$ as $t\to\infty$. We say that the action satisfies \textit{an effective mean ergodic theorem with exponent $\kappa$} if for any $f\in L^2(X)$,
$$\|\beta_{G_t}f-\int_{X}fdm_X\|_2\ll_\kappa \frac{\|f\|_2}{m_G(G_t)^\kappa}.$$
We refer the reader to \cite{GorodnikNevo2015} for a comprehensive survey on effective mean ergodic theorems and their number theoretic applications.

It was shown in  \cite[Theorem 1]{GhoshKelmer17} that if the action of $G$ on $L^2(X)$ satisfies a mean ergodic theorem with exponent $\kappa$ and that the set $G_t=\{g\in G\ |\ \|g\|\leq t\}$ has measure $m_G(G_t)\gg t^b$, then for any $a<2\kappa b$ and any family of shrinking targets $B_t\subseteq X$ with measure $m_X(B_t)\gg t^{-a}$, for a.e. $x\in X$ for all $t\geq t_0(x)$ there is $g\in G_t$ with $g^{-1}x\in B_t$. Taking our space $X=\tilde{G}/\tilde{\G}$ and our shrinking sets $B_t$ as in \lemref{lem:Reduction2} will reduce the problem to establishing a mean ergodic theorem. 

While it is possible to obtain such a mean ergodic theorem in this setting, the exponent $\kappa$ depends on the (strong) spectral gap for the representation of $G$ on $L^2(X)$ which for rank one groups may depend on the lattice $\G$. To remove this dependence we take further advantage of the specific structure of the shrinking sets $B_t$ to give the following refined version.

We can identify $L^2(G/\G)$ with the subspace of $L^2(\tilde{G}/\tilde{\G})$ composed of functions that are invariant under the action of $\R^n$, and let $L^2_{00}(\tilde{G}/\tilde{\G})$ denote its orthogonal complement (that is the set of all function whose average over $\R^n/\Z^n$ is zero). Decomposing $L^2(G/\G)$ further as a direct sum of the space of constant functions and the space, $L^2_0(G/\G)$, of mean zero functions, we get the following decomposition 
\begin{equation}\label{equ:decom}
L^2(\tilde{G}/\tilde{\G})=\C\oplus L_0^2(G/\G)\oplus L^2_{00}(\tilde{G}/\tilde{\G}),
\end{equation}
which is preserved by the left regular $G$-action (under the natural embedding $G\subseteq \tilde{G}$ sending $g\in G$ to $(g,0)\in \tilde{G}$).
 We now reduce the shrinking target problem in \lemref{lem:Reduction2} to bounds for the averaging operators for the latter two representations.
 
\begin{Prop}\label{prop:Reduction3}
Let $G_t\subset G$ be a family of growing sets of measure $m_G(G_t)\gg t^b$ for some $b>0$. Let $\kappa_1>0$ and assume the following 
\begin{enumerate}
\item There is some $\kappa_2>0$ such that for any $f\in L^2(G/\G)$
$$\|\beta_{G_t}f-\int_{G/\G}fdm_G\|_2\ll \frac{\|f\|_2}{m_G(G_t)^{\kappa_2}}.$$
\item For any $\kappa<\kappa_1$ and for any $f\in L^2_{00}(\tilde{G}/\tilde{\G})$
$$\|\beta_{G_t}f\|_2\ll_\kappa \frac{\| f\|_2}{m_G(G_t)^{\kappa}}.$$
\end{enumerate} 
Then for any $a<2b\kappa_1$ and for any family of shrinking targets $A_t\subseteq \bT^n$ with $m(A_t)\gg t^{-a}$, if $B_t\subseteq \tilde{G}/\tilde{\G}$ is as in \lemref{lem:Reduction2}, then
 for a.e. $x\in  \tilde{G}/\tilde{\G}$ and for all sufficiently large $t,$ there is $g\in G_t$ with $g^{-1}x\in B_t$.
\end{Prop}
\begin{proof}
Let $\psi_t\in L^2(\tilde{G}/\tilde{\G})$ denote the indicator function of $B_t$. We note that a Haar measure of $\tilde{G}$ decomposes as $dm_{\tilde{G}}(g,\alpha)=dm_{G}(g)dm(\alpha)$, and by our definition $\tilde{B}_t$ is contained in a single fundamental domain of $\tilde{G}/\tilde{\G}$. Thus 
$$m_{\tilde{G}}(B_t)=m_{\tilde{G}}(\tilde{B}_t) =\int_{\{(g,\alpha)\in\tilde{G}\ |\ g\in \cF_c, \alpha\in \tilde{A}_t\}}dm_G(g)dm(\alpha)=m_G(\cF_c)m(A_t)\asymp_cm(A_t),$$
where $\cF_c=\{g\in G\ | \  g\in \cF, \|g\|\leq c\}$ with $\cF$ the fixed fundamental domain for $G/\G$ as in \lemref{lem:Reduction2}. We want to show that for a.e. $x\in \tilde{G}/\tilde{\G}$ and for all sufficiently large $t,$ there is $g\in G_t$ such that $g^{-1}x\in B_t$. It suffices to show that $\beta_{G_t}(\psi_t)(x)\neq 0$ for all sufficiently large $t$ where $\beta_{G_t}$ denotes the averaging operator \eqref{e:AveOp}. We thus need to show that the set $\cC=\bigcap_{T\geq 0}\bigcup_{t\geq T} \{x\in \tilde{G}/\tilde{\G}\ |\  \beta_{G_t}(\psi_t)(x)=0\}$ has measure zero.
Now we consider the dyadic decomposition 
$$\bigcup_{t\geq T}\{x\in \tilde{G}/\tilde{\G}\ |\ \beta_{G_t}(\psi_t)(x)=0\}= \bigcup_{k\geq \log(T)}\bigcup_{2^k\leq t<2^{k+1}}\{x\in \tilde{G}/\tilde{\G}\ |\ \beta_{G_t}(\psi_t)(x)=0\},$$
and note that, since $G_t$ is increasing and $\psi_t$ is decreasing in $t$, if $\beta_{G_t}(\psi_t)(x)=0$ for some $2^k\leq t< 2^{k+1}$ then %$\beta_{G_t}(\psi_{2^{k+1}})(x)=0$ and hence 
$\beta_{G_{2^k}}(\psi_{2^{k+1}})(x)=0$ so that 
$$\bigcup_{t\geq T}\{x\in \tilde{G}/\tilde{\G}\ |\ \beta_{G_t}(\psi_t)=0\}\subseteq  \bigcup_{k\geq \log(T)}\cC_{2^k,2^{k+1}}$$
where $\cC_{T,t}=\{x\in \tilde{G}/\tilde{\G}\ |\ \beta_{G_T}(\psi_t)(x)=0\}$. We thus need to show that the series $\sum_k m_{\tilde{G}}(\cC_{2^k,2^{k+1}})$ is summable.

We now use our assumptions on the norms of the averaging operators to estimate $m_{\tilde{G}}(\cC_{T,t})$.
Let $\varphi_t\in L^2(G/\G)$ denote the projection of $\psi_t$ and note that $\varphi_t=m(A_t)\chi_{\overline{\cF}_c}$ with $\overline{\cF}_c$ the projection of $\cF_c$ from $G$ to $G/\G$. We also note that $\psi_t-\varphi_t \in L^2_{00}(\tilde{G}/\tilde{\G})$ and 
$\varphi_t-m_{\tilde{G}}(\psi_t)=\varphi_t-m_G(\varphi_t)\in L^2_0(G/\G)$.
Now for any $T,t>1$ we can estimate 
$$\|\beta_{G_T} (\psi_t)-m_{\tilde{G}}(\psi_t)\|_2\leq \|\beta_{G_T}(\psi_t-\varphi_t)\|_2+\|\beta_{G_T}(\varphi_t)-m_G(\varphi_t)\|_2.$$
Now using the bound on the norms of the averaging operators in these spaces we get that for any $\kappa<\kappa_1$
$$\|\beta_{G_T}(\psi_t-\varphi_t)\|_2 \ll_{\kappa} \frac{\|\psi_t-\varphi_t\|_2}{m_G(G_T)^\kappa}\ll_{c} \frac{\sqrt{m(A_t)}}{m_G(G_T)^{\kappa}}.$$  
For the second term, since $\varphi_t=m(A_t)\chi_{\overline{\cF}_c}$ we can bound 
$$\|\beta_{G_T}(\varphi_t)-m_{\tilde{G}}(\psi_t)\|_2=\|\beta_{G_T}(\varphi_t)-m_{G}(\varphi_t)\|_2\ll_c \frac{m(A_t)}{m_{G}(G_T)^{\kappa_2}}.$$
Combining both bounds we get that 
$$\|\beta_{G_T}(\psi_t)-m_{\tilde{G}}(\psi_t)\|_2^2\ll_{\kappa, c} \frac{m(A_t)}{m_G(G_T)^{2\kappa}}+\frac{m(A_t)^2}{m_G(G_T)^{2\kappa_2}}.$$
Since for any $x\in \cC_{T,t}$ we have that $\beta_{G_T}(\psi_t)(x)=0,$  the Chebyshev inequality gives
$$m_{\tilde{G}}(\cC_{T,t})) m(A_t)^2 \ll_c \|\beta_{G_t} \psi_t-m_{\tilde{G}}(\psi_t)\|_2^2,$$
and hence
$$m_{\tilde{G}}(\cC_{T,t})\ll_{\kappa,c} \frac{1}{m(A_t)m_G(G_T)^{2\kappa}}+\frac{1}{m_G(G_T)^{2\kappa_2}}\ll t^{a}T^{-2\kappa b}+T^{-2\kappa_2a}.$$
In particular, assuming that $a<2\kappa_1 b$ we can find $\kappa<\kappa_1$ so that $a<2\kappa b$ for which 
$m_{\tilde{G}}(\cC_{2^k,2^{k+1}})\ll_{\kappa} 2^{k(a-2\kappa b)}+2^{-2k\kappa_2 a}$ is summable, thereby finishing the proof.
\end{proof}

\section{Effective mean ergodic theorems}
In this section we show how the needed effective mean ergodic theorems follow from results on the strong spectral gap. We first recall some general results on such mean ergodic theorems.
\subsection{Relation to operator norms}\label{s:OperatorNorm}
It is useful to think of the averaging operators $\beta_{G_t}$ and the effective mean ergodic theorem in a wider context for general unitary representations.
Given a unitary representation $\pi$ of $G$ on some Hilbert space $\cH$, and a growing family $G_t$ of measurable subsets of $G$ as above, we can consider the averaging operator 
\begin{equation}\label{e:AveOp2}
\beta_{G_t}(\pi)(v)=\frac{1}{m_G(G_t)}\int_{G_t} \pi(g)vdm_G(g).
\end{equation}
This is an operator acting on $\cH$ and we denote by $\|\beta_{G_t}(\pi)\|,$ its operator norm.
Now for the special case where $\pi$ is the representation of $G$ on the space $L^2_0(X)$ of mean zero functions given by $\pi(g)f(x)=f(g^{-1}x)$, the bound $\|\beta_{G_t}(\pi)\|\ll_\kappa m_G(G_t)^{-\kappa}$ is equivalent to an effective mean ergodic theorem with exponent $\kappa$.

One advantage of working with this wider context is that we can reduce the problem to that of irreducible representations.
Explicitly we, record a useful result relating the bound on the operator norm for a representation to that of its irreducible components.
\begin{Lem}\label{lem:reduction}
Let $\pi$ be a unitary representation of $G$ and consider the decomposition $\pi=\int_{Y}^\oplus \pi_{y} d\nu(y)$ as a direct integral of irreducible representations. 
If for some $t>0$ and for $\nu$-a.e. $y\in Y$ the norm of the averaging operator satisfies $\|\beta_{G_t}(\pi_y)\|\leq F(t),$ then $\|\beta_{G_t}(\pi)\|\leq F(t)$.
\end{Lem}
\begin{proof}
From our assumption for $\nu$-a.e. $y\in Y,$ for all $v_y\in \cH_y$, we have that  $\|\beta_{G_t}(\pi_y) v_y\|\leq F(t)\|v_y\|$.
Since for any $v\in \cH$ we have $\|v\|^2=\int_Y \|v_y\|^2d\nu(y)$ and $(\beta_{G_t}(\pi)v)_y=\beta_{G_t}(\pi_y) v_y$, we get that 
$$\|\beta_{G_t}(\pi)v\|^2=\int_Y\| \beta_{G_t}(\pi_y) v_y\|^2d\nu(y)\leq F(t)^2\|v\|^2,$$
so that the operator norm satisfies that $\|\beta_{G_t}(\pi)\|\leq F(t)$ as claimed.
\end{proof}

Now, for a unitary representation $\pi$, there is a close relation between bounds for operator norms of $\beta_{G_t}(\pi)$ and the strong spectral gap for $\pi$. The strong spectral gap is closely related to decay of matrix coefficients and is controlled by the parameter $p(\pi)\in [2,\infty)$ as introduced in the introduction.
% namely the infimum of all $p\geq 2$ for which there exists a dense set of vectors $v\in \pi$ such that the matrix coefficient $g\mapsto \langle\pi(g)v, v\rangle$ lies in $L^p(G)$. Recall also that $\pi$ is tempered when $p(\pi)=2$.
%(see \cite{KelmerSarnak2009} for more background). 
It follows from the work of Gorodnik and Nevo \cite{GorodnikNevo2010} (see also \cite[p. 78]{GorodnikNevo2015} and \cite[p. 306]{Nevo1998}) that for any unitary representation $\pi$ of a semi-simple Lie group $G$, if we let $l$ be the smallest even integer satisfying $l\geq p(\pi)/2$ and let $\kappa_1=\frac{1}{2l}$  (when $\pi$ is tempered we can take $\kappa_1=1/2$) then
for any $\kappa<\kappa_1$  
\begin{equation}\label{e:SG2MET}
\|\beta_{G_t}(\pi)\|\ll_\kappa m_G(G_t)^{-\kappa}.
\end{equation}

\begin{rem}
We note that when the growing sets $G_t$ are bi-$K$-invariant, it is expected that in fact one can take $\kappa_1=\frac{1}{p(\pi)}$ without the parity conditions.
In particular, when $G$ has real rank one this follows from estimates on decay of matrix coefficients of spherical functions. While this result is well known to experts, as we could not find a reference in the litrature, for the sake of completeness we include a proof below.  
\end{rem}
\begin{Prop}\label{p:effp1}
Keep the notation as above and assume that $Q$ is of signature $(n-1,1)$. Let $\{G_t\}_{t>1}$ be the family of growing norm balls defined in \eqref{e:G_T}. Then for any $\kappa<\frac{1}{p(\pi)}$ we have 
$$\|\beta_{G_t}(\pi)\|\ll_{\kappa} m_G(G_t)^{-\kappa}.$$
%In particular for $\pi$ the representation of $G$ on $L^2_{00}(\tilde{G}/\tilde{\G})$ this holds for any $\kappa<\frac{1}{2(n-2)}$.
\end{Prop}
\begin{proof}
Note that if $Q,Q'$ are two forms of signature $(n-1,1)$ then there is $h\in \SL_n(\R)$ and $\lambda\neq 0$ with $Q'(v)=\lambda Q(vh)$ and hence conjugating by $h$ gives an isomorphism $\varphi_h: G'\to G$ with $G'=\SO_{Q'}^+(\R)$.
Now, any unitary representation $\pi'$ of $G'$ is of the form $\pi'=\pi\circ \varphi_h$ and it is clear that in this case $p(\pi)=p(\pi')$. We also have that the corresponding norm balls satisfy  
$G_{t/c}\subseteq G'_t\subseteq G_{ct}$ for some $c>1$ from which it follows that it is enough to prove the result for a single form, and we may take $Q=Q_0$ given in \eqref{equ:baseform}.

Next, we recall that when $G$ is of rank one, the decay of matrix coefficients can be given explicitly in terms of the $KA^+K$ decomposition with $A^+=\{a_t\ |\ t\geq 0 \}$ the positive chamber in a Cartan subgroup, and $K$ the corresponding maximal compact subgroup, as in \secref{subsec:coordinate}. More precisely, for a unitary representation $\pi$ of $G$ on a Hilbert space $\cH$ not weakly containing the trivial representation, let $\alpha(\pi)\in (0,\frac{n-2}{2}]$ be the smallest number satisfying that for any positive $\alpha<\alpha(\pi)$ and for any $K$-finite vector $v\in\cH$, we have for any $t>0$%In particular, we can also characterize the spectral gap of a unitary representation $\pi$ on a Hilbert space  $\cH$ by the parameter $\alpha(\pi)$ which is the supremum over all $\theta>0$ such that any $K$-finite vector $v\in \cH$ satisfies 
\begin{equation}\label{equ:decay}
\langle \pi(a_t)v,v\rangle \ll_{\alpha} \textrm{dim$\langle \pi(K)v\rangle$}e^{-\alpha t}\|v\|^2.
\end{equation} 
Using  \eqref{equ:Haar} for signature $(n-1,1)$ we have that,  up to scaling, the Haar measure of $G$ in the coordinates $g=ka_tk'$ is given by 
$$dm_G(g)=\sinh(t)^{n-2}dtdkdk',$$
which gives the relation that 
\begin{equation}\label{equ:reductionrelation}
\alpha(\pi)=\frac{n-2}{p(\pi)}.
\end{equation}

Now for any test vector $v$, let $v^K=\int_K \pi(k)vdk$. Note that $\|v^K\|\leq \|v\|$ and since $G_t$ is bi-$K$-invariant then 
$$\beta_{G_t}(\pi)v=\beta_{G_t}(\pi)v^K.$$
Hence to calculate the operator norm we just need to estimate $\|\beta_{G_t}(\pi)v\|$ for $v$ a spherical vector.
Now for $v$ a norm one spherical vector, using the $KA^+K$ decomposition, the estimate \eqref{e:mG_T}, the description of $G_t$ \eqref{equ:gtp1}, and the fact that $\pi$ is unitary, we have for $t>T_0$
\begin{align*}
\|\beta_{G_t}(\pi)v\|^2&\leq \frac{1}{m_G(G_t)^2}\int_{G_t}\int_{G_t}\left|\langle \pi(g_1)v,\pi(g_2)v\rangle\right| dm_G(g_1)dm_G(g_2)\\
&\asymp \frac{1}{t^{2(n-2)}}\int_0^{\log(t)}\int_0^{\log(t)}\int_K\left|\langle \pi(a_{t_1}k a_{t_2})v,v\rangle\right| \sinh(t_1)^{n-2} \sinh(t_2)^{n-2}dkdt_1dt_2.
\end{align*}
Further decomposing $k=mk_\theta m'$ with $m,m'\in M$ and $k_\theta$ as in \eqref{e:ktheta}, and noting that $m, m'$ commute with $a_t$ and using the Haar measure decomposition in \eqref{e:HaarMKM},  we get that the second line of the above equation is given, up to a constant, by 
\begin{align*}
\asymp_n\frac{1}{t^{2(n-2)}}\int_0^{\log(t)}\int_0^{\log(t)}\int_0^\pi \left|\langle \pi(a_{t_1}k_\theta a_{t_2})v,v\rangle\right| \sin(\theta)^{n-3}\sinh(t_1)^{n-2} \sinh(t_2)^{n-2}d\theta dt_1dt_2.
\end{align*}
Now we use the $KA^+K$ decomposition to write $a_{t_1}k_{\theta}a_{t_2}=ka_tk'$ for some $k,k'\in K$ and $a_t\in A^+$, and use the decay of matrix coefficients to estimate the above integral. For any positive $\alpha<\alpha(\pi)$, using \eqref{equ:decay}, we can estimate matrix coefficients of a spherical norm one vector by
$$|\langle \pi(a_{t_1}k_{\theta}a_{t_2})v,v\rangle|=|\langle \pi(ka_t k')v,v\rangle|\ll_{\alpha} e^{-\alpha t}\asymp \cosh(t)^{-\alpha}.$$
We thus need to estimate the term $\cosh t$ in terms of the coordinates $t_1$ and $t_2$. By comparing the $(n,n)^{\rm{th}}$ entry of both matrices we see that
\begin{equation}\label{equ:relaaka}
\cosh(t)=\cosh(t_1)\cosh(t_2)+\cos(\theta)\sinh(t_1)\sinh(t_2).
\end{equation}
We can rearrange this, noting that $\cosh(t_1)\cosh(t_2)-\sinh(t_1)\sinh(t_2)=\cosh(t_1-t_2)$ to get that 
$$\cosh(t)=2\cos^2(\theta/2)\sinh(t_1)\sinh(t_2)+\cosh(t_1-t_2)\geq 2\cos^2(\theta/2)\sinh(t_1)\sinh(t_2).$$
Using this bound and the estimate on decay of matrix coefficients we can estimate
\begin{align*}
\|\beta_{G_t}(\pi)v\|^2\ll_{\alpha}  \frac{1}{t^{2(n-2)}}\int_0^{\log(t)}\int_0^{\log(t)}\int_0^\pi \left(\cos^2(\theta/2)\sinh(t_1)\sinh(t_2)\right)^{-\alpha}  \\
\sin(\theta)^{n-3}\sinh(t_1)^{n-2} \sinh(t_2)^{n-2}d\theta dt_1dt_2.
\end{align*}
For the innermost integral over $\theta$, note that since $\alpha<\alpha(\pi)=\frac{n-2}{p(\pi)}\leq \frac{n-2}{2}$ the integral 
$$\int_0^\pi \cos(\theta/2)^{-2\alpha} \sin(\theta)^{n-3}d\theta$$
converges. We can thus estimate
\begin{align*}
\|\beta_{G_t}(\pi)v\|^2\ll_{\alpha} \frac{1}{t^{2(n-2)}}\left(\int_0^{\log(t)} \sinh(t_1)^{n-2-\alpha}dt_1\right)^2\ll_{\alpha} t^{-2\alpha}.
\end{align*}
Since this holds for any norm one spherical vector, we get that for any positive $\alpha<\alpha(\pi)=\frac{n-2}{p(\pi)}$,
$$\|\beta_{G_t}(\pi)\|\ll t^{-\alpha}\asymp m_G(G_t)^{-\frac{\alpha}{n-2}}.$$
In particular, for any $\kappa<\frac{1}{p(\pi)}$ we can take $\alpha=(n-2)\kappa<\alpha(\pi)$ to conclude the proof. 
\end{proof}
%\begin{remark}\label{rmk:decrel}
%Let $G=\SO_Q^+(\R)$ with $Q$ a general quadratic form of signature $(n-1,1)$ and let $\pi$ be a unitary representation of $G$. Then there is $h\in \SL_n(\R)$ such that $h^{-1}Gh=\SO_{Q_0}^+(\R)$. The Haar measure of $G$ is the push-forward of the Haar measure of $\SO_{Q_0}^+(\R)$. Using the above Haar measure description of $\SO_{Q_0}^+(\R)$ and the relation $\|a_t\|\asymp e^t$ it is not difficult to see that if a function $F$ on $G$ satisfies $|F(g)|\ll \|g\|^{-\alpha}$, then $F\in L^p(G)$ for all $p>(n-2)/\alpha$.
%\end{remark}

\subsection{Groups with property $(T)$}
For a connected semi-simple Lie group $G$ with finite center, define 
$$p(G):=\sup\{p(\pi)\ |\ \textrm{$\pi$ is a nontrivial irreducible unitary representation of $G$}\}.$$
We note that $G$ has property $(T)$ if and only if $p(G)<\infty$. Thus for groups with property $(T)$ we can bound $p(\pi)$ from above by $p(G)$ for any unitary representation $\pi$ of $G$ not containing the trivial representation.
Effective bounds for $p(G)$ were obtained for all semi-simple Lie groups with property $(T)$ in  \cite{Li1995,Oh98}, implying in particular the following 
\begin{Prop}\label{prop:spectralgap}
For $G=\SO^+_Q(\R)$ with $Q$ a form of signature $(p,q)$ we have 
\[p(G)\leq \left\{
  \begin{array}{ll}
    p+q-2 & p+q\geq 7,\; q\geq 2 \textrm{ and } (p,q)\not\in \{(5,2),(4,3), (6,3)\}\\
    2(p-1) & (p,q)\in\{(5,2),(4,3), (6,3)\} \\
   6 & (p,q)\in \{(4,2), (3,3)\}\\
   4 & (p,q)=(3,2).\\  
  \end{array}
\right.
\]
\end{Prop}
\begin{proof}
In all the above cases the group $G$ has property $(T)$. When $(p,q)\notin\{(5,2),(4,3),(6,3)\}$ the parameter $p(G)$ was explicitly computed in \cite{Li1995} yielding the cases except $(p,q)\in \{(5,2),(4,3),(6,3)\}$. For the remaining cases this bound follows from the upper bound on $p(G)$ proved by Oh \cite{Oh98}. 
\end{proof}
Combining these bounds for the strong spectral gap gives the following result on an effective mean ergodic theorem for these groups.
\begin{Thm}\label{t:METPropT}
Let $G=\SO^+_Q(\R)$ with $Q$ a form in $n\geq 5$ variables of signature $(p,q)$ with $p\geq q> 1$. Let 
$$\kappa_1=\left\lbrace\begin{array}{cc}
\frac{1}{n} &  n\equiv 0\pmod{4}\\
\frac{1}{n-1} &  n\equiv 1\pmod{4}\\
\frac{1}{n-2} &  n\equiv 2\pmod{4}\\
\frac{1}{n+1} & n\equiv 3\pmod{4}\\
\end{array}\right.$$
except for the case of signatures $(3,3)$, $(4,2)$ and $(6,3)$ for which we let $\kappa_1=1/8$, $1/8$ and $1/12$ respectively.
Then for any $\kappa<\kappa_1,$ for any unitary representation $\pi$ of $G$ without nontrivial $G$-invariant vectors, and for any growing family $G_t$ with finite positive measure we have that  
$$\|\beta_{G_t}(\pi)\|\ll_\kappa m_G(G_t)^{-\kappa}.$$
\end{Thm}

\subsection{Signature $(n-1,1)$}
For a form of signature $(n-1,1)$ the group $G=\SO^+_Q(\R)$ does not have property $(T)$ so there is no uniform bound for the strong spectral gap in general. However, in order to apply \propref{prop:Reduction3}, the representations we are interested in are the representations of $G$ on the two function spaces $L^2_0(G/\G)$ and $L^2_{00}(\tilde{G}/\tilde{\G})$  described in the decomposition \eqref{equ:decom}.
%is the restriction to $G$ of a representation of the semi-direct product $\tilde{G}=G\ltimes \R^n$. That is, the group $\tilde{G}$ acts on the space $L^2(\tilde{G}/\tilde{\G})$ and preserves the decomposition \eqref{equ:decom}.
%$$L^2(\tilde{G}/\tilde{\G})=L^2(G/\G)\oplus L^2_{00}(\tilde{G}/\tilde{\G}).$$

Using  \thmref{t:rank1a} (whose proof is postponed to the next section) we have that $p(\tilde\Gamma)_0=\max\{2,n-2\}$. Hence using \propref{p:effp1} assumption (2) of \propref{prop:Reduction3} holds with $\kappa_1=\frac{1}{2}$ for $n=3,4$ and $\kappa_1=\frac{1}{n-2}$ for $n\geq 4$.

\begin{rem}\label{EMETRank1}
For $Q$ of signature $(n-1,1)$ we don't have a uniform bound for the strong spectral gap for $L^2(G/\G)$ (unless $\G$ is a congruence group). However, in this case the spectral gap is equivalent to the strong spectral gap and the discreteness of the Laplacian spectrum implies that there is some bound for the spectral gap (which may depend on $\G$). We thus get an effective mean ergodic theorem of the form as assumption (1) in \propref{prop:Reduction3} with some exponent $\kappa_2>0$ that may depend on $\G$.
\end{rem}

\subsection{Signature $(2,2)$} \label{s:22}
When $Q$ is of signature $(2,2)$, the group $\SO_Q(\R)$ is locally isomorphic to $\SL_2(\R)\times \SL_2(\R)$ and does not possess property $(T)$.
To see this local isomorphism more clearly, it will be convenient to work with the determinant form 
\begin{equation}\label{e:Q1}
Q_1(a,b,c,d)=ad-bc=\det(M)
\end{equation}
 when identifying $\R^4=\mathrm{Mat}_2(\R)$ and writing $M=\left(\begin{smallmatrix} a & b\\ c& d\end{smallmatrix}\right)$. 
Consider the action of $\SL_2(\R)\times \SL_2(\R)$ on $\R^4=\mathrm{Mat}_2(\R)$ with $(g_1,g_2)\in \SL_2(\R)\times \SL_2(\R)$ sending $M\in \R^4$ to $g_1Mg_2^*$. This action is clearly linear and preserves $Q_1$, and thus induces a homomorphism from $\SL_2(\R)\times \SL_2(\R)$ to $\SO_{Q_1}(\R)$. In fact, let $G=\SO_{Q_1}^+(\R)$ be the identity component of $\SO_{Q_1}(\R)$, then this action induces a double covering 
\begin{equation}\label{e:iota}
\iota: \SL_2(\R)\times \SL_2(\R)\to G
\end{equation}
with the kernel $\ker(\iota)=\{\pm(I_2,I_2)\}$. Thus any irreducible unitary representation of $G$ is of the form
$\pi(\iota(g_1,g_2))=\pi_{1}(g_1)\otimes \pi_2(g_2)$ where each $\pi_i$ is an irreducible unitary representation of $\SL_2(\R)$ such that $\pi_1\otimes \pi_2$ is trivial on $\ker(\iota)$. By a slight abuse of notation we will write in this case 
$\pi=\pi_1\otimes\pi_2$ and we will view $\pi$ as a representation of both $\SL_2(\R)\times \SL_2(\R)$ and $G$. 

Nevertheless, \thmref{t:SG}  (whose proof we postpone to the next section) implies that the representation of $G$ on $L^2_{00}(\tilde{G}/\tilde{\G})$  is tempered.
This has the following immediate corollary, allowing us to obtain an optimal exponent for the effective mean ergodic theorem for functions in $L^2_{00}(\tilde{G}/\tilde{\G})$ for any growing family of sets $G_t$. 
\begin{Cor}\label{cor:sig22}
Keep the notation as above. For any growing family of sets $G_t$ in $G$ we have for all $\kappa<1/2$, for any $f\in L^2_{00}(\tilde{G}/\tilde{\G})$, and for all $t>1$,
$$\|\beta_{G_t}f\|_2 \ll_\kappa \frac{\|f\|_2}{m_G(G_t)^\kappa}.$$
\end{Cor}

The situation for functions in $L^2(G/\G)$ is more complicated. 
Here the lattice $\G$ is not necessarily an irreducible lattice 
(for example for $Q=Q_1$ above we can identify $G$ with $\SL_2(\R)\times\SL_2(\R)$ and then $\G=\SL_2(\Z)\times\SL_2(\Z)$). In particular, the representation of $G$ on the space 
$L^2(G/\G)$ might not have a strong spectral gap.  Hence, in order to get an effective rate for the mean ergodic theorem we need to make sure the growing norm balls are well balanced.

Before defining the well balanced norm balls, we analyze the norm balls $G_t \subset G$ defined in \eqref{e:G_T} using the coordinates from the double cover $\iota: H\times H\to G=\SO_{Q_1}^+(\R)$ with $H=\SL_2(\R)$. Fix a polar decomposition $H=\SO_2(\R)A^+\SO_2(\R)$ with 
$$A^+=\left\{a_t=\left(\begin{smallmatrix}
e^{t/2} & 0\\
0 & e^{-t/2}\end{smallmatrix}\right)\ |\ t>0\right\},$$
and for any $g\in H$ we denote by $t(g)>0$ the uniquely determined positive number in the decomposition $g=k_1a_{t(g)}k_2$ with $k_1,k_2\in \SO_2(\R)$ and $a_{t(g)}\in A^+$. We note that in these coordinates the Haar measure of $H$ is given, up to a scalar, by
$$dm_H(k_1a_tk_2)=\sinh(t) dk_1dtdk_2$$
with $dk$ the probability Haar measure on $\SO_2(\R)$. Moreover, for $g=k_1a_{t(g)}k_2$ as above we have $g^{-1}=k_2^{-1}a_{t(g)}^{-1}k_1^{-1}=k_2^{-1}\omega a_{t(g)}\omega^{-1}k_1^{-1}$ implying that $t(g)=t(g^{-1})$. Here $\omega=\left(\begin{smallmatrix} 0 & -1\\
1 & 0\end{smallmatrix}\right)\in \SO_2(\R)$.

Let $\|\cdot\|$ be the Euclidean norm on $\R^4$. Note that under the identification $\R^4= \mathrm{Mat}_2(\R)$, $\|\cdot\|$ is the Hilbert-Schmidt norm on $\mathrm{Mat}_2(\R)$ which is bi-$\SO_2(\R)$-invariant. First, for any $t_1,t_2>0$ and any $M=\left(\begin{smallmatrix} a & b\\ c& d\end{smallmatrix}\right)\in \mathrm{Mat}_2(\R)$ we have
$$a_{t_1}M a^*_{t_2}=\begin{pmatrix} ae^{(t_1+t_2)/2} & be^{(t_1-t_2)/2}\\ ce^{(t_2-t_1)/2}& de^{-(t_1+t_2)/2}\end{pmatrix},$$
implying that $\|\iota(a_{t_1},a_{t_2})\|_{\rm{op}}=e^{(t_1+t_2)/2}$, where the operator norm is attained when taking $M=\left(\begin{smallmatrix}
1 & 0\\
0 &0 \end{smallmatrix}\right)$. Now for any $\iota(g_1,g_2)\in G$ and for any $M=\left(\begin{smallmatrix} a & b\\ c& d\end{smallmatrix}\right)\in \mathrm{Mat}_2(\R)$ writing $g_i=k_i'a_{t_i}k_i$ with $k_i,k_i'\in \SO_2(\R)$ and $a_{t_i}\in A^+$ for $i=1,2$ we have
$$\|g_1Mg_2^*\|=\|a_{t_1}k_1Mk_2^*a_{t_2}\|\leq e^{(t_1+t_2)/2}\|k_1Mk_2^*\|=e^{(t_1+t_2)/2}\|M\|,$$
implying that $\|\iota(g_1,g_2)\|_{\rm{op}}\leq e^{(t_1+t_2)/2}$. On the other hand taking $M$ such that $k_1Mk_2^*=\left(\begin{smallmatrix}
1 & 0\\
0 &0 \end{smallmatrix}\right)$ in the above equation we get 
$$\|g_1Mg_2^*\|=e^{(t_1+t_2)/2}\|k_1Mk_2^*\|=e^{(t_1+t_2)/2}\|M\|,$$
implying that $\|\iota(g_1,g_2)\|_{\rm{op}}\geq e^{(t_1+t_2)/2}$. Hence we have $\|\iota(g_1,g_2)\|_{\rm{op}}= e^{(t_1+t_2)/2}=\|\iota(g_1,g_2)^{-1}\|_{\rm{op}}$ (since $t(g)=t(g^{-1})$ for any $g\in H$) implying that
$$G_T=\{\iota(g_1,g_2)\in G\ |\ t(g_1)+t(g_2)\leq 2\log T \}.$$
Recall that by \eqref{e:mG_T} we have $m_G(G_T)\asymp T^2\log T$. Denoting by
 $$H_T:=\{g\in H\ |\ t(g)\leq 2\log T\},$$ 
 we see that the projection to each factor $G_T\cap \left(H\times \{I\}\right)=H_T$ has measure $\asymp T^2$, and hence,  the growing norm balls are balanced but are not well balanced in the sense of \cite[Definition 3.17]{GorodnikNevo2010}.
We thus need to replace the norm balls with slightly smaller well balanced norm balls given by 
$$G^{\rm wb}_T=\{\iota(g_1,g_2)\in G_T\ |\ \max\{t_1(g),t_2(g)\}\leq \log(T)\}=\iota\left(H_{\sqrt{T}}\times H_{\sqrt{T}}\right)$$
having measure $m_G( G^{\rm wb}_T)=m_H(H_{\sqrt{T}})^2\asymp T^2$. 

For $G=\SO_Q^+(\R)$ with $Q$ a different quadratic form of signature $(2,2)$, we fix a conjugation isomorphism $\varphi_Q : \SO_{Q_1}^+(\R)\to G$ and define $G^{\rm wb}_T:=\varphi_Q(\iota(H_{\sqrt{T}}\times H_{\sqrt{T}}))$. We note that since $G^{\rm wb}_T\subset G_T$ when $G=\SO_{Q_1}^+(\R)$, for general $G=\SO_Q^+(\R)$ we have $G^{\rm wb}_T\subset G_{cT}$ for some constant $c>1$ depending only on $Q$.

Now, for these well balanced balls we can show the following.
\begin{Thm}\label{t:NB}
For  $G=\SO^+(2,2)$ and $\{G^{\rm wb}_T\}_{T>1}$ the well balanced norm balls defined above,
there is some $\kappa>0$ $($that may depend on $\G$$)$ such that
 $\|\beta_{G^{\rm wb}_T}(f)\|_2\ll \frac{\|f\|_2}{m_G(G^{\rm wb}_T)^{\kappa}}$, for any $f\in L^2_0(G/\G)$. 
\end{Thm}
\begin{proof}
In view of the reduction arguments in the proof of \propref{p:effp1}, it suffices to prove this theorem for the case when $G=\SO_{Q_1}^+(\R)$ with $Q_1$ given in \eqref{e:Q1}. From \lemref{lem:reduction} it is enough to show that 
$\|\beta_{G^{\rm wb}_T}(\pi)\|\ll \frac{1}{m_G(G^{\rm wb}_T)^{\kappa}}$ for any irreducible component $\pi=\pi_1\otimes \pi_2$ of $L^2_0(G/\G)$ with each $\pi_i$ irreducible representation of $\SL_2(\R)$.

Recall that any non trivial irreducible unitary representation of $\SL_2(\R)$ is infinitesimally equivalent to one of the following: the spherical and non-spherical principal series representations, the discrete series representations, the two mock discrete series representations and the complementary series representations, see e.g. \cite[Chapter \rom{6}]{LangSL2R} for more details on the description of the unitary dual of $\SL_2(\R)$.   We note that among these irreducible representations the only non-tempered representations are the complementary series representations, and following the parameterization in \cite{LangSL2R} up to infinitesimal equivalence they can be parameterized by the interval $(0,1)$. We thus denote them by $\sigma_s$ with $s\in (0,1)$, and by examining the decay rate of matrix coefficients (see e.g. \cite[p. 216]{HoweTan1992}) and using the relation \eqref{equ:reductionrelation} we see that $\sigma_s$ has spectral gap $p(\sigma_s)=\frac{2}{1-s}$. We also denote the trivial representation of $\SL_2(\R)$ by $\sigma_1$.

Now, the averaging operator for the balanced norm balls takes the form 
\begin{align*}
\beta_{G^{\rm wb}_T}(\pi)v&=\frac{1}{m_G( G^{\rm wb}_T)}\int_{H_{\sqrt{T}}\times H_{\sqrt{T}}}\pi_1(h_1)\otimes \pi_2(h_2) vdm_H(h_1)dm_H(h_2)\\
&=\frac{1}{m_H(H_{\sqrt{T}})}\int_{H_{\sqrt{T}}}\pi_1(h_1)\left(\frac{1}{m_H(H_{\sqrt{T}})}\int_{H_{\sqrt{T}}} \pi_2(h_2) vdm_H(h_2)\right)dm_H(h_1)\\
&= \beta_{H_{\sqrt{T}}}(\pi_1)\beta_{H_{\sqrt{T}}}(\pi_2)v.
\end{align*}
Now, if one of the two representations, say $\pi_1$, is tempered then for any $\kappa_1<1/2$ we have 
$\|\beta_{H_{\sqrt{T}}}(\pi_1)\|\ll_{\kappa_1} m_H(H_{\sqrt{T}})^{-\kappa_1}$. Hence 
\begin{align*}
\|\beta_{ G^{\rm wb}_T}(\pi)v\|&=\| \beta_{H_{\sqrt{T}}}(\pi_1)\beta_{H_{\sqrt{T}}}(\pi_2)v\|\\
&\ll_{k_1} m_H(H_{\sqrt{T}})^{-\kappa_1}\|\beta_{H_{\sqrt{T}}}(\pi_2)v\|\leq m_H(H_{\sqrt{T}})^{-\kappa_1}\|v\|,
\end{align*}
where we used the trivial bound $\|\beta_{H_{\sqrt{T}}}(\pi_2)\|\leq 1$ for $\pi_2$. Since $m_H(H_{\sqrt{T}})=m_G(G^{\rm wb}_T)^{1/2}$, this proves the claim in this case with $\kappa=\frac{\kappa_1}{2}<1/4$.

It remains to treat the case where both representations are non-tempered. The discreteness of the Laplacian spectrum 
on $L^2(G/\G)$ implies that there is some constant $s_0\in (0,2)$ (which may depend on $\G$) such that any irreducible representation $\pi_1\otimes\pi_2$ weakly contained in $L^2_0(G/\G)$ satisfies that if both $\pi_i=\sigma_{s_i}$ with $s_i\in (0,1]$ are either complementary series or the trivial representation, then $s_1+s_2\leq s_0$. In particular, for such $\pi$ at least one of $\pi_i$ is nontrivial, say, $\pi_2=\sigma_{s_2}$ is always nontrivial.

For the norm balls $H_{\sqrt{T}}$ in $H=\SL_2(\R)$ and the complementary series representation $\sigma_s$ (with $p(\sigma_s)=\frac{2}{1-s}$) by \propref{p:effp1} (noting that $H=\SL_2(\R)$ is locally isomorphic to $\SO(2,1)$) we have
$\|\beta_{H_{\sqrt{T}}}(\sigma_{s_i})\|\ll_{\kappa_i} m_H(H_{\sqrt{T}})^{-\kappa_i}$ for any $\kappa_i<\frac{1-s_i}{2}$ (we use the trivial bound $\|\beta_{H_{\sqrt{T}}}(\sigma_{s_1})\|\leq 1$ if $\sigma_{s_1}=\sigma_1$ is the trivial representation). Now let $0<\kappa<\frac{1}{2}-\frac{s_0}{4}$ and we can take $\kappa_i<\frac{1-s_i}{2}$ (we take $\kappa_1=0$ if $s_1=1$) such that $(\kappa_1+\kappa_2)/2=\kappa$ to get that
\begin{align*}
\|\beta_{G^{\rm wb}_T}(\pi)v\|&=\| \beta_{H_{\sqrt{T}}}(\pi_1)\beta_{H_{\sqrt{T}}}(\pi_2)v\| \ll_{\kappa_1} m_H(H_{\sqrt{T}})^{-{\kappa_1}}\|\beta_{H_{\sqrt{T}}}(\pi_2)v\| \\
&\ll_{\kappa_2} m_H(H_{\sqrt{T}})^{-(\kappa_1+\kappa_2)}=m_G(G^{\rm wb}_T)^{-(\kappa_1+\kappa_2)/2}= m_G(G^{\rm wb}_T)^{-\kappa}
\end{align*}
as claimed. 
\end{proof}

\subsection{Proof of main results}
Collecting together \thmref{t:METPropT}, \propref{p:effp1} and \corref{cor:sig22} gives the proof of \thmref{t:EMET}, where for signature $(2,2)$ we use the well balanced normed balls. The proof of \thmref{thm:genericshifts} then follows as described in \secref{sec:reduction}. More precisely, by \eqref{e:mG_T}, we have that $m_{G}(G_t) \gg t^{q(p-1)}$, while in signature $(2,2)$ the same estimate holds for $m_G(G_t^{\rm wb})$. So by \propref{prop:Reduction3}, \lemref{lem:Reduction1} and \lemref{lem:Reduction2}, we have that  the conclusion of \thmref{thm:genericshifts} holds for $\kappa_0 = 2\kappa_1q(p-1)$. Note that the first condition of \propref{prop:Reduction3} trivially holds when $G$ has property $(T)$, while in the remaining cases it follows from the discreteness of the spectrum of the Laplace operator (see \rmkref{EMETRank1} and \thmref{t:NB}). 
%The estimation of $\kappa_1$ is carried out in \thmref{t:METPropT} for forms in $5$ or more variables of signature $(p, q)$ with $p \geq q > 1$. For forms of signature $(n-1,1)$ the exponent $\kappa_1$ follows from the result on the spectral gap \thmref{thm:spectralgapp1} together with  \thmref{thm:effp1}  relating the spectral gap to the mean ergodic theorem. 
Again, for forms of signature $(2,2)$ we use \propref{prop:Reduction3} with the well balanced norm balls $G_t^{\rm{wb}}$ instead of $G_t$, noting that $G_t^{\rm{wb}}\subset G_{ct}$ for some $c>1$.

\section{Upper and lower bounds for $p(\tilde{\Gamma})_0$}
Let $G=\SO_Q^+(\R), \tilde{G}=G\ltimes \R^n$, $\G\leq \SO_Q^+(\Z)$ a finite-index subgroup and $\tilde{\G}=\G\ltimes \Z^n$ be as above.  
Throughout this section we denote by $\pi$ the representation of $G$ on $L^2(\tilde{G}/\tilde{\Gamma})$ which is unitary with respect to the inner product 
$\langle f_1,f_2\rangle=\int_{\tilde{G}/\tilde{\G}}f_1(x)\overline{f_2(x)}dm_{\tilde{G}}(x)$. Our goal is to give upper and lower bounds for $p(\tilde{\G})_0=p(\pi)$.
For both upper and lower bounds we will need the following general construction of functions in $L^2_{00}(\tilde{G}/\tilde{\G})$.
\begin{Lem}\label{l:testfunctions}
Given a bounded and compactly supported function $\vf$ on $G$ and a nontrivial character $\lambda$ of $\bT^n=\R^n/\Z^n$, let $f=f_{\vf,\lambda}$ be given by
\begin{equation}\label{e:flambda}
f(g,v)=\sum_{\gamma\in \G} \vf(g\gamma)\lambda(v\gamma).\end{equation}
Then $f\in L^2_{00}(\tilde{G}/\tilde{\Gamma})$, and for any $F\in  L^2(\tilde{G}/\tilde{\Gamma})$
\begin{equation}\label{e:unfolding}
\langle f,F\rangle= \int_{G} \vf(g)\int_{\bT^n }\lambda(v)\overline{F(g,v)}dm(v)dm_G(g).
\end{equation}
In particular, the family of functions 
$f_{\vf,\lambda}$ with $\lambda$ a nontrivial character and $\varphi$ compactly supported, span a dense subspace of $L^2_{00}(\tilde{G}/\tilde{\Gamma})$.
\end{Lem}
\begin{proof}
The identity $f(g\gamma,v\gamma+u)=f(g,v)$ for any $(\g,u)\in \G\ltimes \Z^n$ is straight forward, and since $\varphi$ is bounded and compactly supported (so that its support can be covered by finitely many fundamental domains for $G/\G$), we have that $f$ is bounded, implying that $f\in L^2(\tilde{G}/\tilde{\G})$. Moreover, since $\int_{\bT^n}\lambda(v\g)dm(v)=0$ for any $\g\in \G$, then $\int_{\bT^n}f(g,v)dm(v)=0$. Hence $f\in L^2_{00}(\tilde{G}/\tilde{\G})$. 
Next to show \eqref{e:unfolding} let $\cF$ be a fundamental domain for $G/\G$ and $\cP$ a fundamental parallelogram for $\bT^n$, and write
%\begin{eqnarray*}
%f((g,v)(\gamma_0, u_0))&=&=\sum_{\g\in \G}\vf_0(g\gamma_0\gamma)\lambda(v\gamma_0\gamma+u\gamma)\\
%&=& \sum_{\g\in \G}\vf_0(g\gamma)\lambda(v\gamma)=f(g,v)\\
%\end{eqnarray*}
% and if $\cF$ is a fundamental domain for $G/\Gamma$ and $\cP$ a fundamental parallelogram for $\bT^n$ then 
 \begin{align*}
 \langle f,F\rangle &= \int_{\tilde{G}/\tilde{\G}}f(g,v)\overline{F(g,v)}dm_G(g)dm(v)\\
 &=\sum_{\g\in \G}\int_{\cF}\int_{\cP}\vf(g\gamma)\lambda(v\gamma)\overline{F(g,v)}dm(v)dm_G(g)\\
 &=\sum_{\g\in \G}\int_{\cF\gamma} \vf(g)\int_{\cP \gamma}\lambda(v)\overline{F(g,v)}dm(v)dm_G(g).
\end{align*}
 Since  for any $g$ the function $v\mapsto \lambda(v)F(g,v)$ is a function on $\bT^n$ and $\cP\gamma$ is also a fundamental domain for $\R^n/\Z^n$ we may replace $\cP\gamma$ by $\cP$ above to get 
  \begin{align*}
 \langle f,F\rangle &=\sum_{\g\in \G}\int_{\cF\gamma} \vf(g)\int_{\cP}\lambda(v)\overline{F(g,v)}dm(v)dm_G(g)\\
 &=\int_{G} \vf(g)\int_{\cP }\lambda(v)\overline{F(g,v)}dm(v)dm_G(g)
\end{align*}
as claimed. %In particular, taking $F=f$ which is bounded shows that  $f\in L^2_{00}(\tilde{G}/\tilde{\G})$. \comm{[Doesn't $f$ is bounded automatically imply $f\in L^2(\tilde{G}/\tilde{\G})$ so we don't need to apply the above integral formula.]}

Finally for the density argument, for any $F\in L^2_{00}(\tilde{G}/\tilde{\G})$  let 
$$F_\lambda(g)=\int_{\bT^n}\overline{F(g,v)}\lambda(v)dm(v),$$
 and note that $F$ is orthogonal to $f_{\vf,\lambda}$ for all compactly supported $\vf$ implies that 
$$\int_G \vf(g)F_\lambda(g)dm_G(g)=0,$$ 
for all compactly supported $\vf$, and hence $F_\lambda=0$.
Now, if $F\in L^2_{00}(\tilde{G}/\tilde{\G})$ is orthogonal to all functions $f_{\vf,\lambda}$  with nontrivial $\lambda$ and $\vf$ compactly supported, then it satisfies that 
$F_\lambda(g)=0$ for all characters $\lambda$ and hence $F=0$.
\end{proof}
\subsection{Lower bounds}
To get lower bounds for $p(\tilde\Gamma)_0$ we first prove the following result, giving an upper bound for the critical exponent in the mean ergodic theorem, by presenting an explicit $f\in L^2_{00}(\tilde{G}/\tilde{\G})$ for which $\|\beta_{G_t}f\|_2$ is large. 
We first give a large family of test functions for which we have an explicit estimate for the norm of the averaging operator.
\begin{Prop}\label{prop:lower}
Let $\lambda$ be a nontrivial character of $\bT^n=\R^n/\Z^n$ and let $G^\lambda\leq G$ denote its stabilizer
%\footnote{Here the notation $G_{\lambda}$ looks similar to the norm ball notation. Since we use $\lambda$ exclusively in this paper as a character on $\R^n/\Z^n$, we hope it does not cause confusion.}
%\comm{[Do you think this footnote suffice? Alternatively, we can change the norm balls to $B_t$.]}
$$G^\lambda=\{g\in G: \lambda(vg)=\lambda(v),\ \forall\ v\in \R^n\}.$$
Then for any  $f=f_{\vf,\lambda}$ as in \eqref{e:flambda} with $\vf$ non-negative and compactly supported, there is a constant $c>0$ depending only on the support of $\vf$ such that for all $t>c$,
$$\frac{\left|\int_{G}\vf dm_G\right|^2}{m_G(G_t)^2}\sum_{\g\in \G\cap G^\lambda}m_G(\g G_{t/c}\cap G_{t/c})\leq \|\beta_{G_t}f\|_2^2\leq \frac{\left|\int_{G}\vf dm_G\right|^2}{m_G(G_t)^2}\sum_{\g\in \G\cap G^\lambda}m_G(\g G_{ct}\cap G_{ct}).$$
\end{Prop}
%\begin{remark}
%T
%\end{remark}
\begin{proof}
%For simplicity of notation, for $g\in G$ and $v\in \R^n$ we abbreviate $dm_G(g)$ and $dm(v)$ to $dg$ and $dv$ respectively.
Starting from 
\begin{align*}
\|\beta_{G_t}f\|_2^2&=\frac{1}{m_G(G_t)^2}\int_{G_t}\int_{G_t} \langle \pi(h_1)f,\pi(h_2)f\rangle dm_G(h_1)dm_G(h_2)\\
&= \frac{1}{m_G(G_t)^2}\int_{G_t}\int_{G_t} \langle f,\pi(h_1^{-1}h_2)f\rangle dm_G(h_1)dm_G(h_2)
\end{align*}
and using \eqref{e:unfolding} with $F(g,v)=\pi(h_1^{-1}h_2)f(g,v)=f(h_2^{-1}h_1g,v)$ we get that 
\begin{align*}
\|\beta_{G_t}f\|_2^2 &= \frac{1}{m_G(G_t)^2}\int_{G_t}\int_{G_t} \int_G \vf(g)\int_{\bT^n} \overline{f(h_2^{-1}h_1 g,v)}\lambda(v)dm(v)dm_G(g) dm_G(h_1)dm_G(h_2).
\end{align*}
Further expanding $f(g,v)=\sum_{\g\in \G} \vf(g\g)\lambda(v\g)$ and noting that $\int_{\bT^n}\lambda(v)\overline{\lambda(v\g)}dm(v)=0$ unless $\g\in G^\lambda$  we get that 
\begin{align*}
\|\beta_{G_t}f\|_2^2% &=& \sum_{\g\in \G}\frac{1}{m_G(G_t)^2}\int_{G_t}\int_{G_t} \int_G \vf(g) \vf(h_2^{-1}h_1 g\g)\int_{\bT^n}\lambda(v\g)\lambda(v)dm(v)dm_G(g) dm_G(h_1)dm_G(h_2)\\
&=  \sum_{\g\in \G\cap G^\lambda}\frac{1}{m_G(G_t)^2}\int_{G_t}\int_{G_t} \int_G \vf(g) \overline{\vf(h_2^{-1}h_1 g\g)}dm_G(g) dm_G(h_1)dm_G(h_2)\\
&=  \sum_{\g\in \G\cap G^\lambda}\frac{1}{m_G(G_t)^2}\int_{G_t}\int_{G_t} \int_G \vf(h_1^{-1}g)\overline{ \vf(h_2^{-1} g\g)}dm_G(g) dm_G(h_1)dm_G(h_2).
\end{align*}
Let $\chi_{G_t}$ denote the indicator function of $G_t$ then making a change of variables $h_1\mapsto gh_1$ and $h_2\mapsto  g\g h_2$ gives 
\begin{align*}
\|\beta_{G_t}f\|_2^2&= \frac{1}{m_G(G_t)^2}\sum_{\g\in \G\cap G^\lambda}\int_{G}\int_{G}\int_{G}\chi_{G_t}(h_1)\chi_{G_t}(h_2) \vf(h_1^{-1}g)\overline{\vf(h_2^{-1}g\g)}dm_G(h_1)dm_G(h_2)dm_G(g)\\
&= \frac{1}{m_G(G_t)^2}\sum_{\g\in \G\cap G^\lambda}\int_{G}\int_{G}\int_{G}\chi_{G_t}(gh_1)\chi_{G_t}(g\g h_2) \vf(h_1^{-1}) \overline{\vf(h_2^{-1})}dm_G(h_1)dm_G(h_2)dm_G(g).
\end{align*}
Since $\vf$ is compactly supported there is some $c>0$ such that $\max\{\|h\|, \|h^{-1}\|\}\leq c$ for all $h\in G$ with $h^{-1}$ in the support of $\vf$. 
%we have that $\|h^{-1}\|\leq c/\tilde{C}$ implying that (cf. \eqref{equ:comnorm}) $\max\{\|h\|, \|h^{-1}\|\}\leq c$. 
This then further implies that for such $h$,
%\comm{[the first inequality is not exactly true: it's equivalent to $G_{t/c}\subset G_{t}h^{-1}$, or $G_{t/c}h\subset G_t$, but we only have $\|h^{-1}\|\leq c$; we may also need to use $\|h\|\asymp_n\|h^{-1}\|$]}
$$\chi_{G_{t/c}}(g)\leq \chi_{G_t}(gh)\leq \chi_{G_{ct}}(g),\quad \forall\ g\in G.$$
Since we assume $\vf$ is non-negative we get the lower bound 
\begin{align*}
\|\beta_{G_t}f\|_2^2&\geq \ \frac{1}{m_G(G_t)^2}\sum_{\g\in \G\cap G^\lambda}\int_{G}\int_{G}\int_{G}\chi_{G_{t/c}}(g)\chi_{G_{t/c}}( g\g) \vf (h_1^{-1})\vf(h_2^{-1})dm_G(h_1)dm_G(h_2)dm_G(g)\\
%&= \frac{1}{m_G(G_t)^2}\sum_{\g\in \G\cap H_\lambda}\int_{G}\chi_{G_{t/c}}(g)\chi_{G_{t/c}}(\g g)dm_G(g)\left(\int_{G}\vf(h)dm_G(h)\right)^2\\
&= \frac{\left(\int_{G}\vf(h)dh\right)^2}{m_G(G_t)^2}\sum_{\g\in \G\cap G^\lambda}m_{G}( G_{t/c}\cap G_{t/c}\g).
\end{align*}
The upper bound follows from the same argument. 
%\begin{align*}
%\|\beta_{G_t}f\|_2^2&\leq \ \frac{1}{m_G(G_t)^2}\sum_{\g\in \G\cap H_\lambda}\int_{G}\int_{G}\int_{G}\chi_{G_{ct}}(g)\chi_{G_{ct}}(\g g) \vf (h_1^{-1})\overline{\vf(h_2^{-1})}dm_G(h_1)dm_G(h_2)dm_G(g)\\
%&= \frac{1}{m_G(G_t)^2}\sum_{\g\in \G\cap H_\lambda}\int_{G}\chi_{G_{ct}}(g)\chi_{G_{ct}}(\g g)dm_G(g)\left|\int_{G}\vf(h)dm_G(h)\right|^2\\
%&= \frac{\left|\int_{G}\vf(h)dh\right|^2}{m_G(G_t)^2}\sum_{\g\in \G\cap H_\lambda}m_{G}(\g G_{ct}\cap G_{ct}).
%\end{align*}
%The lower bound follows from the same argument recalling that $\vf$ is real valued and non negative.
\end{proof}

In order to use this formula to estimate the norm of the averaging operator we need a good estimate for 
$m_{G}(G_{t}\cap G_{t}\g)$. However, since we are only interested in a lower bound, the following simple estimate will do.
\begin{Lem}\label{lem:cbound}
Let $A\leq G$ be a Cartan subgroup and let $K\leq G$ be a maximal compact subgroup such that $G=KA^+K$. Let $C=\sup\{\|k\|\ |\ k\in K\}$. Then we have that $m_{G}(G_{t}\cap G_{t}\g)\gg 1$ uniformly for all $t>2C$ and for all $\g\in \G$ with $\|\g\|\leq t^2/(4C^2)$.
\end{Lem}
\begin{proof}
Let  $\g\in \G$ with $\|\g\|\leq t^2/(4C^2)$. Decomposing $\g=kak'$ with $k,k'\in K$ and $a\in A$, we can find some $a'\in A$ satisfying $\|a'\|\leq t/(2C)$ and $\|a' a\|= \frac{\|a\|}{\|a'\|}\leq t/(2C)$ (if $\|a\|\leq t/(2C)$ then we can take $a'$ to be the identity element, and if $\|a\|>t/(2C)$ then we can take $a'$ such that $\|a' a\|= \frac{\|a\|}{\|a'\|}= t/(2C)$ and in both cases we have $\|a'\|\leq t/(2C)$). Let $g_0=a'k^{-1}$ so that $\|g_0\|=\|a'k^{-1}\|\leq C\|a'\|\leq t/2$ and $\| g_0\g\|=\|a'ak'\|\leq C\|a' a\|\leq t/2$. Hence any $g\in G_{2}g_0$  will satisfy that $\|g\|\leq t$ and $\|g\g \|\leq t$, implying that 
$m_{G}( G_{t}\cap G_{t}\g)\geq m_G(G_2)$.
\end{proof}

Note that $m_G(G_{ct})\asymp m_G(G_t)$ and, since $\G\cap G^\lambda$ is a lattice in $G^\lambda$,  by \cite[Theorem 6.4]{GorodnikNevo2010} we have that for large $t$
$$\#\{\g\in \G\cap G^\lambda\ |\ \|\g\|\leq t\}\asymp m_{G^\lambda}(G^\lambda \cap G_t).$$
Using these estimates leads to the following.
\begin{Cor}
Let $G^\lambda\subseteq G$ be the stabilizer of a nontrivial character $\lambda$ of $\bT^n=\R^n/\Z^n$, then there is $f\in L^2_{00}(\tilde{G}/\tilde{\G})$ such that 
$$\|\beta_{G_t}f\|_2^2\gg_{f} \frac{m_{G^\lambda}(G^\lambda \cap G_{t^2})}{m_G(G_t)^2}.$$
In particular, since for $G\cong \SO^+(p,q)$ with $p>2$ and $p\geq q$ we have a character $\lambda$ with stabilizer $G^\lambda$ locally isomorphic to $\SO^+(p-1,q)$ %\comm{[Need to check: is the stabilizer really the identity component, or could it be $\SO(p-1,q)$]} 
we get a function $f=f_{\vf,\lambda}$ for which 
$$
\|\beta_{G_t}f\|_2^2\gg \frac{m_{G^\lambda}(G^\lambda\cap G_{t^2})}{m_G(G_t)^2}\gg\left\lbrace\begin{array}{cc}  t^{2(p-2)q-2(p-1)q}
\gg m_{G}(G_t)^{-\frac{2}{p-1}}& p>q\\
t^{2(p-1)^2-2(p-1)p}\left(\log t\right)^{-2}\gg_\e m_G(G_t)^{-\frac{2}{p}-\e} & p=q.
\end{array}\right.
$$
\end{Cor}

Combining this result and the relation between $p(\pi)$ and the critical exponent in the mean ergodic theorem we get the following.
\begin{proof}[Proof of  \thmref{t:LowerBound}]
First note that if the representation $\pi$ of $G$ on $L^2_{00}(\tilde{G}/\tilde{\Gamma})$ was tempered then 
$\|\beta_{G_t}f\|_2\leq \frac{\|f\|_2}{m_G(G_t)^{\kappa}}$
for any $\kappa<\tfrac{1}{2}$ and any $f\in L_{00}^2(\tilde{G}/\tilde{\Gamma})$. However, if $p>q$ and $p>3$ then $\frac{1}{p-1}<\tfrac{1}{2}$ while for if $p=q\geq 3$ then we can find $\e$ sufficiently small so that $\frac{1}{p}+\tfrac{\e}{2}<\tfrac{1}{2}$.
This implies that the representation is not tempered except possibly for signatures $(3,2),(3,1),(2,2)$ and $(2,1)$. %(For signatures $(2,1)$, $(3,1)$ and $(2,2)$ we know that $\pi$ is tempered following from \cite{}, 

Next for the case of signature $(n-1,1)$ with $n\geq 4$, since for all $\kappa<\frac{1}{p(\tilde{\G})_0}$ we have the bound  $\|\beta_{G_t}f\|_2\leq \frac{\|f\|_2}{m_G(G_t)^{\kappa}}$, the example above implies that  $p(\tilde\Gamma)_0\geq n-2$.
For signature $(n-q,q)$ with $1<q<\frac{n}{2}$ let $j\in\{0,1,2,3\}$ with $j\equiv n-q-2 \pmod{4}$. Then our example implies that $p(\tilde{\G}_0)> n-q-2-j$. Indeed if $p(\tilde{\G}_0)\leq n-q-2-j$ then $\|\beta_{G_t}f\|_2\leq \frac{\|f\|_2}{m_G(G_t)^{\kappa}}$  for any $\kappa<\frac{1}{n-q-2-j}$ which contradicts our counter example. Finally, for signature $(q,q)$ the same argument shows that $p(\pi)>n-q-1-j$  with $j\in\{0,1,2,3\}$ satisfying $j\equiv n-q-1\pmod{4}$.
\end{proof}

\subsection{Upper bounds.}
In order to get upper bound for $p(\pi)=p(\tilde\G)_0$ we need bounds for matrix coefficients of a dense set of test functions.
 We start with the following general estimate.
\begin{Lem}\label{l:MatrixCoefficient}
For $\vf,\vf'$ compactly supported functions on $G$ and $\lambda,\lambda'$ two nontrivial characters of $\bT^n=\R^n/\Z^n$ let $f=f_{\vf,\lambda}$ and $f'=f_{\vf',\lambda'}$ be as in \eqref{e:flambda}.
Then for any $h\in G$ we have that  $\langle f', \pi(h)f\rangle=0$ unless $(\g_0)_*\lambda=\lambda$ for some $\gamma_0\in \Gamma$, in which case
\begin{align*}
|\langle f',\pi(h)f\rangle|\leq \|\vf\|_\infty \|\vf'\|_\infty  \sum_{\g\in \G\cap G^\lambda} m_G( hG_c\cap G_c \g_0\g),
\end{align*}
where $\gamma_{*}\lambda(v):=\lambda(v\gamma)$ for any $v\in \R^n$,   $G^\lambda$ is the stabilizer of $\lambda$ in $G$ given as in \propref{prop:lower}, and $c\geq 1$ satisfies that $\vf, \vf'$ are supported on $G_c$.
\end{Lem}
\begin{proof}
Using \eqref{e:unfolding} with $F=\pi(h)f$ and expanding $f(g,v)$ we get that 
\begin{align*}
\langle f',\pi(h)f\rangle&=\int_G \vf'(g)\int_{\bT^n} \overline{f(h^{-1}g,v)}\lambda'(v)dm(v)dm_G(g)\\
&=\sum_{\g\in \G}\int_G \vf'(g)  \overline{\vf(h^{-1}g\g)}\int_{\bT^n}\overline{\lambda(v\gamma)}\lambda'(v)dm(v)dm_G(g)\\
&=\mathop{\sum_{\g\in \G}}_{\gamma_*\lambda=\lambda'}\int_G \vf'(g)  \overline{\vf(h^{-1}g\g)}dm_G(g).
\end{align*}
%where $\gamma_{*}\lambda(v):=\lambda(v\gamma)$ for any $v\in \R^n$.
The sum is empty unless $(\gamma_0)_*\lambda=\lambda'$ for some $\g_0\in \G$ in which case making a change of variables $\g\mapsto \g_0\g$ we get that  
\begin{align*}
\langle f',\pi(h)f\rangle&=
\mathop{\sum_{\g\in \G\cap G^{\lambda}}}\int_G \vf(hg)  \overline{\vf(g\g_0\g)}dm_G(g).
\end{align*}
The result now follows by bounding $|\vf| \leq \|\vf\|_\infty\chi_{G_c}$  and $|\vf'|\leq \|\vf'\|_\infty\chi_{G_c}$.
\end{proof}

In order to use this estimate we need to control the size of $m_G( hG_c\cap G_c \g_0\g)$, which can be done as follows when $Q$ is of signature $(n-1,1)$.
\begin{Lem}\label{l:hGGh}
Let $G=\SO^+_Q(\R)$ with $Q$ of signature $(n-1,1)$. Then for any $h_1,h_2\in G$, %and let $h_1,h_2\in G$ of the form $h_1=k_1a_{t_1}k_1'$ and $h_2=k_2a_{t_2}k_2'$ correspond to the $KA^+K$ decomposition.
%Then
$m_G(h_1G_c\cap G_ch_2)=0$ unless $\|h_1\|\asymp_c\|h_2\|$ in which case
$m_G(h_1G_c\cap G_c h_2)\ll_c \|h_1\|^{-(n-2)}$.
\end{Lem}
%\comm{[We didn't really define the Cartan decomposition for a general form, so here I rephased this lemma in terms of the norm and also added remark 4.7 for the next proof.]}
\begin{proof}
Arguing similarly as in the beginning of the proof of \propref{p:effp1} we may assume $Q=Q_0$ as in \eqref{equ:baseform}. In particular, in this case the norm balls are bi-$K$-invariant. Writing $h_1=k_1a_{t_1}k_1'$ and $h_2=k_2a_{t_2}k_2'$ and using the invariance of the Haar measure we then have $m_G(h_1G_c\cap G_ch_2)=m_G(a_{t_1}G_c\cap G_ca_{t_2})$, and that
\begin{align}\label{e:AKK}
\nonumber m_G(a_{t_1}G_c \cap G_ca_{t_2})&=\int_G \chi_{G_c}(a_{t_1}g)\chi_{G_c}(ga_{t_2})dm_G(g)\\
&=\int_{K}\int_{0}^\infty \int_{K} \chi_{G_c}(a_{t_1}k_1a_t)\chi_{G_c}(a_tk_2a_{t_2})(\sinh(t))^{n-2}dtdk_1dk_2\\
\nonumber &=\int_{0}^\infty \left(\int_{K} \chi_{G_c}(a_{t_1}ka_t)dk\right)\left(\int_K\chi_{G_c}(a_tk a_{t_2})dk\right)(\sinh(t))^{n-1}dt.
\end{align}
We thus need to estimate the two inner integrals. Further decompose $k=m_1k_\theta m_2$ with $k_\theta$ defined as in \eqref{e:ktheta} and $m_1,m_2$ commuting with $a_t$, to get that 
$$\int_{K} \chi_{G_c}(a_{t_1}ka_t)dk=\int_0^\pi \chi_{G_c}(a_{t_1}k_\theta a_t)\sin(\theta)^{n-3}d\theta.$$
We can now use  \eqref{equ:relaaka} to write $a_{t_1}k_\theta a_{t}=k_1 a_\tau k_2$ with 
$$\cosh(\tau)=2\cos^2(\theta/2)\sinh(t_1)\sinh(t)+\cosh(t_1-t).$$
Using that $\| k_1 a_\tau k_2\|\asymp \cosh(\tau)$ the condition that  $a_{t_1}k_\theta a_t\in G_c$ implies that
$$2\cos^2(\theta/2)\sinh(t_1)\sinh(t)+\cosh(t_1-t)\ll c.$$ In particular $t=t_1+O_c(1)$ and $\cos^2(\theta/2)\ll_c  e^{-2t_1}$.
Consequently, the product $$\left(\int_{K} \chi_{G_c}(a_{t_1}ka_t)dk\right)\left(\int_K\chi_{G_c}(a_tk a_{t_2})dk\right)=0,$$ unless $t_1=t+O_c(1)$ and $t_2=t+O_c(1)$.
Hence  $m_G(a_{t_1}G_c \cap G_ca_{t_2})=0$ unless $t_1=t_2+O_c(1)$, or equivalently, $\|h_1\|\asymp_c\|h_2\|$. 
 In this case, the condition $\cos^2(\theta/2)\ll_c  e^{-2t_1}$ implies $\theta=\pi+O_c(e^{-t_1})$ so $\sin(\theta)\ll e^{-t_1}$ and 
 $$\int_{K} \chi_{G_c}(a_{t_1}ka_t)dk=\int_0^\pi \chi_{G_c}(a_{t_1}k_\theta a_t)\sin(\theta)^{n-3}d\theta\ll e^{-t_1(n-2)}\asymp \|h_1\|^{-(n-2)}.$$
Similarly, we also have that $\int_K\chi_{G_c}(a_tk a_{t_2})dk\ll_c \|h_2\|^{-(n-2)}$. Plugging these estimates, together with the estimates $t=t_1+O_c(1)=t_2+O_c(1)$ back in \eqref{e:AKK} completes the proof.
\end{proof}

We can now calculate the spectral gap $p(\tilde{\G})_0$ for any finite-index subgroup $\G\leq \SO^+_Q(\Z)$ with $Q$ of signature $(n-1,1)$.
\begin{proof}[Proof of \thmref{t:rank1a}]
Let $G=\SO_Q^+(\R)$ with $Q$ of signature $(n-1,1)$, let  $\G\leq \SO^+_Q(\Z)$ denote a finite-index subgroup, and let $\tilde\G=\G\ltimes \Z^n$ be as above. We have already shown in \thmref{t:LowerBound} that $p(\tilde{\G})_0\geq n-2$;
it thus remains to show that $p(\tilde{\G})_0\leq n-2$. To do this we need to find a dense set of vectors in $L^2_{00}(\tilde{G}/\tilde{\G})$ for which the corresponding matrix coefficients are in $L^p(G)$ for all $p>n-2$.
We will use the set of vectors spanned by the functions $f_{\vf,\lambda}$ given as in \eqref{e:flambda} with $\vf$ compactly supported and $\lambda$ a nontrivial character of $\bT^n$. Note that by \lemref{l:testfunctions} this is indeed a dense set in $L^2_{00}(\tilde{G}/\tilde{\G})$. From linearity it is enough to show that for any $p>n-2$ the functions $g\mapsto \langle f_{\vf',\lambda'}, \pi(g)f_{\vf,\lambda}\rangle$ are all in $L^p(G)$.
First, by  \lemref{l:MatrixCoefficient} it is enough to consider the case when $(\gamma_0)_*\lambda=\lambda'$ for some $\g_0\in \G$ in which case we can estimate
\begin{align*}
|\langle f_{\vf',\lambda'}, \pi(g)f_{\vf,\lambda}\rangle|\leq \|\vf\|_\infty \|\vf'\|_\infty  \sum_{\g\in \G\cap G^\lambda} m_G( gG_c\cap G_c \g_0\g),
\end{align*}
where $G^\lambda$ is the stabilizer of $\lambda$ in $G$, and $c$ is any constant so that $\vf$ and $\vf'$ are supported on $G_c$. Next, %for any $g\in G$ let $t(g)$ correspond to the decomposition $g=k_1a_{t(g)}k_2$ then 
by \lemref{l:hGGh} the term $m_G( gG_c\cap G_c \g_0\g)$ vanishes unless $\|g\|\asymp_c\|\g_0\g\|$ in which case it is bounded by $O_c(\|g\|^{-(n-2)})$. This leads to the bound 
\begin{align*}
|\langle f_{\vf',\lambda'}, \pi(g)f_{\vf,\lambda}\rangle|  \ll_{\vf,\vf',\lambda, \lambda'}   \|g\|^{-(n-2)}\# \{\g\in \G\cap G^\lambda\ |\  \|\g\|\ll_{c,\g_0} \|g\|\}.
\end{align*}
Here for the estimate we also used that $\|\g_0\g\|\asymp_{\g_0}\|\g\|$ and absorbed the dependence on $c, \g_0$ into the dependence on $\vf, \vf',\lambda, \lambda'$.
%where we used that $t(\g_0\g)=t(\g)+O_{\g_0}(1)$.
Now using a classification of the stabilizer group $G^{\lambda}$ (see \propref{p:stabilizers} below) we have that $G^\lambda$ is either compact (in which case $\#(\G\cap G^\lambda)$  is uniformly bounded) or it is a semi-direct product of a maximal unipotent and a compact group, in which case 
$$\# \{\g\in \G\cap G^\lambda\ |\ \|\g\|\ll_{c,\g_0} \|g\|\} \ll \vol\left(\{v\in \R^{n-2}\ |\ \|v\|^2 \ll_{c,\gamma_0} \|g\|\}\right)\ll_{c,\gamma_0} \|g\|^{\frac{n-2}{2}},$$ or it is isomorphic to a copy of $\SO^+(n-2,1)$ inside $G$, in which case 
$$\# \{\g\in \G\cap G^\lambda\ |\ \|\g\|\ll_{c,\g_0} \|g\|\} \ll_{c,\g_0} \|g\|^{n-3}.$$ In particular, in all cases we get that
\begin{align*}
|\langle f_{\vf',\lambda'}, \pi(g)f_{\vf,\lambda}\rangle |&\ll_{\vf,\vf',\lambda,\lambda'}  \|g\|^{-1},
\end{align*}
Finally, since there is $h\in \SL_n(\R)$ such that $h^{-1}Gh=\SO_{Q_0}^+(\R)$ and the Haar measure of $G$ is the push-forward of the Haar measure of $\SO_{Q_0}^+(\R)$. Using  \eqref{equ:Haar}  and the relation $\|a_t\|\asymp e^t$ it is not difficult to see that if a function $F$ on $G$ satisfies $|F(g)|\ll \|g\|^{-1}$, then $F\in L^p(G)$ for all $p>(n-2)$ thus concluding the proof.
\end{proof}

\section{Bounds for $p(\tilde{G})_0$}
In this section we take a closer look at the parameter $p(\tilde{G})_0$ for $G=\SO_Q^+(\R)$ without property $(T)$, explicitly when $Q$ has signature $(2,2)$ or $(n-1,1)$.
\subsection{ Signature $(2,2)$}
%It remains to establish the results on the spectral gap for the cases where $G$ does not have property $(T)$.
%Let $G$ be one of the groups $\SO^+(n-1,1)$ or $\SO^+(2,2)$ and let $\tilde{G}=G\ltimes \R^n$.  We want to give a uniform bound for the strong spectral gap for representations $\pi$ that are restrictions to $G$ of a representation $\tilde{\pi}$ of $\tilde{G}$ having no non-trivial $\R^n$-invariant vectors. We recall that when $G=\SO^+(2,1)\cong \PSL_2(\R)$ by \cite[Theorem 7.3.9]{Zimmer1984}  any such representation is tempered so $p(\pi)=2$ in this case. We will bootstrap this result to give uniform bounds for the strong spectral gap for  $\SO^+(n-1,1)$ and $\SO^+(2,2)$.

%\subsection{Signature $(2,2)$}
Since for different forms of the same signature the corresponding stabilizers are conjugate it is enough to show this for the specific form $Q=Q_1$ given in \eqref{e:Q1}.  In this case we can identify the stabilizer $G=\SO_Q^+(\R)$ with $\SL_2(\R)\times \SL_2(\R)$ where the action of 
$(g_1,g_2)\in \SL_2(\R)\times \SL_2(\R)$ on $\R^4=\mathrm{Mat}_2(\R)$ is given by  $M\mapsto g_1Mg_2^*$.
As we noted in \secref{s:22}, the irreducible representations of $G$ are all of the form $\pi_1\otimes \pi_2$ with $\pi_1,\pi_2$ irreducible representations of $\SL_2(\R)$. With this identification in mind we have the following.

\begin{proof}[Proof of \thmref{t:SG}]
Let  $G=\SO_{Q_1}^+(\R)$ and let $\iota: \SL_2(\R)\times\SL_2(\R)\to G$ be the homomorphism defined in \eqref{e:iota}. Consider the two subgroups $\tilde{G}_1,\tilde{G}_2$ of $\tilde{G}$ given by 
$$\tilde{G}_i=\{(\iota(g_1,g_2),v)\in \tilde{G}\ |\ g_i=I_2\},$$
and let $G_1,G_2\cong \SL_2(\R)$ be the corresponding two subgroups of $G$.
Note that each of the groups  $\tilde{G}_i$ is naturally isomorphic to $\SL_2(\R)\ltimes \R^4$ where the action of $\SL_2(\R)$ on $\R^4=\mathrm{Mat}_2(\R)$ is given by matrix multiplication (one acting on the left and the other acting by the transpose on the right). In particular, for both cases
the only $\SL_2(\R)$-invariant vector is the zero vector.

Now, let $\tilde{\pi}$ be a representation of $\tilde{G}$ with no nontrivial $\R^4$-invariant vectors. Then $\tilde{\pi}|_{\tilde{G}_i}$ is a representation of $\SL_2(\R)\ltimes \R^4$ with no nontrivial $\R^4$-invariant vectors and hence $\tilde{\pi}|_{G_i}$ is tempered (see \cite[Theorem 7.3.9]{Zimmer1984}).  
Now to see that $\pi=\tilde{\pi}|_G$ is tempered, it is enough to show that any irreducible representation weakly contained in $\pi$ is tempered. But any such irreducible representation is of the form $\pi_1\otimes \pi_2$ with $\pi_1,\pi_2$ irreducible representations of $G_1,G_2\cong \SL_2(\R)$. Since the restriction of $\pi$ to each of the factors is tempered we must have that both $\pi_1,\pi_2$ are tempered, and hence $\pi_1\otimes \pi_2$ is tempered. Since this holds for any irreducible representation weakly contained in $\pi$ then $\pi$ is tempered as claimed.
\end{proof}

\subsection{Signature $(n-1,1)$}
For signature $(n-1,1)$ we use a different strategy, using an induction argument. However, in order to execute the induction argument we need to take a closer look at the proof of  \cite[Theorem 7.3.9]{Zimmer1984}, and in particular Mackey's characterization of representations of semi-direct products (see \cite[Theorem 7.3.1]{Zimmer1984})

\begin{thm*}[Mackey]\label{T:Mackey}
Let $G$ be a group acting on $\R^n$ and let $\tilde{G}=G\ltimes \R^n$. For any irreducible unitary representation $\tilde{\pi}$ of $\tilde{G}$  there is a unitary character $\lambda$ of $\R^n$, and an irreducible unitary representation $\sigma$ of $\tilde{G}^\lambda$, 
the stabilizer of $\lambda$ in $\tilde{G}$, such that 
\begin{enumerate}
\item $\tilde\pi=\ind_{\tilde{G}^\lambda}^{\tilde{G}}(\sigma)$,
\item $\sigma |_{\R^n}=(\dim \sigma)\lambda$,
\item $\tilde \pi |_{\R^n}\cong L^2(\tilde{G}/\tilde{G}^\lambda,\cH)$ for some Hilbert space $\cH$ with respect to the measure on $\tilde{G}/\tilde{G}^\lambda$ coming from Haar measure on $\tilde{G}$, where the action of $\R^n$ is given by
$$\left(\tilde\pi(1,v)f\right)(g)=\lambda(vg)f(g).$$
\end{enumerate}
\end{thm*}
\begin{rem}
The group $\tilde{G}$ acts on $\R^n\leq \tilde{G}$ by conjugation and induces an action on the group of unitary characters $\hat{\R}^n$. Since the action of $\R^n$ is trivial this action factors through the group $G$, which acts on characters by $g\cdot \lambda(v)=\lambda(v g)$. In particular, the stabilizer is $\tilde{G}^\lambda=G^\lambda\ltimes \R^n$ with 
$$G^\lambda=\{g\in G\ |\ \lambda(vg)=\lambda(v),\; \forall v\in \R^n\}.$$
Consequently we can identify the quotients  $\tilde{G}/\tilde{G}^\lambda=G/G^\lambda$.
\end{rem}

We use this characterization for the case of $G=\SO_Q^+(\R)$, for $Q$ of signature $(n-1,1)$. To further understand this characterization of irreducible representations of $\tilde{G}$, we take a closer look at the structure of the stabilizers $G^\lambda$ for characters $\lambda$ of $\R^n$.
Any unitary character $\lambda$ of $\R^n$ is of the form $\lambda(v)=e^{i v\cdot  \alpha}$ for some vector $\alpha\in \R^n$, and with this identification we have that %\comm{[we used the notation $G_{\lambda}$ for this stabilizer before.]}
\begin{equation}\label{equ:stab}
G^\lambda=\{g\in G\ |\  \alpha g^* =\alpha\},
\end{equation}
is the transpose of the stabilizer of $\alpha$ in $G^*$ (the transpose of $G$) under the right multiplication action of $G^*$ on $\R^n$. So the first step in understanding the representation $\tilde\pi$ is to characterize the different stabilizers. 

\begin{rem}\label{rmk:trans}
Note that for $G=\SO^+_Q(\R)$, we have its transpose $G^*=\SO_{Q^*}^+(\R)$ with $Q^*$ a different form of the same signature. Explicitly, if $Q(v)=vJv^*$ for some symmetric matrix $J$ with $\det(J)\neq 0$ then $Q^*(v)=vJ^{-1}v^*$ has stabilizer $\SO_{Q^*}(\R)=(\SO_Q(\R))^*$. We thus need to understand the structure of stabilizers in $G^*$.
\end{rem}

First for $\alpha=0$ the character $\lambda$ is trivial, $\tilde{G}^\lambda=\tilde{G}$ and $\tilde{\pi}=\sigma$.
In this case, the restriction $\tilde{\pi}|_{\R^n}=\sigma|_{\R^n}$ is trivial, so this case does not occur when $\pi$ has no nontrivial $\R^n$-invariant vectors.
Next for $\alpha\neq 0$ we note that, up to conjugation in $G^*$, the stabilizer of $\alpha\neq 0$ in $G^*=\SO^+_{Q^*}(\R)$ only depends on the sign of $Q^*(\alpha)$. 
The following proposition summarizes the different possible stabilizers, we omit the proof which is a simple calculation.

%\comm{[shall we move this proposition earlier? we have already used information about these stabilizers before.]}
\begin{Prop}\label{p:stabilizers}
Let $\lambda(v)=e^{iv\cdot\alpha}$ with $\alpha\neq 0$ and $G^\lambda$ as above.
If $Q^*(\alpha)<0$ then $G^\lambda$ is compact; if $Q^*(\lambda)=0$, then $G^\lambda$ is conjugate to the semi-direct product of a maximal unipotent subgroup and a compact group; and if $Q^*(\alpha)>0$, then  $G^\lambda$ is  a copy of $\SO^{+}_{Q'}(\R)$ sitting inside $G$ with $Q'$ a form of signature $(n-2,1)$ (given by the restriction of $Q$ to $V_\lambda$).
\end{Prop}

For cases where the stabilizer $G^\lambda$ is amenable we also have that $\tilde{G}^\lambda=G^\lambda\ltimes \R^n$ is amenable. Hence in these cases $\sigma$ is weakly contained in the regular representation of $\tilde{G}^\lambda$ and hence $\tilde\pi$ is weakly contained in the regular representation of $\tilde{G}$. But then any irreducible component of $\pi=\tilde{\pi}|_{G}$ is weakly contained in the regular representation of $G$ and hence is tempered. So in these cases we have that $p(\pi)= 2$. We note that when $n=3$ the stabilizer $G^\lambda$ is always  amenable and all representations are tempered, however, when $n>3$ this is no longer the case when $Q^*(\alpha)>0$.

To handle the cases where $G^\lambda$ is not amenable, let $V_\lambda=\ker(\lambda)=\{v\in \R^n\ |\ \lambda(v)=1\}$ and identify the semi-direct product 
$G^\lambda \ltimes V_\lambda$ as a subgroup of $\tilde{G}=G\ltimes \R^n$. We then show the following. 
\begin{Lem}\label{l:Vlambda}
Keep the notation as above and assume that $Q^*(\alpha)\neq 0$. Then
the representation $\tilde{\pi}$ has no nontrivial $V_\lambda$-invariant vectors.
\end{Lem}
\begin{proof} 
It is enough to show that the restriction $\tilde{\pi}|_{\R^n}$ has no nontrivial $V_\lambda$-invariant vectors, and from the characterization of  $\tilde \pi |_{\R^n}\cong L^2(\tilde{G}/\tilde{G}^\lambda,\cH)$ it is enough to show that for any  $f\in L^2(\tilde{G}/\tilde{G}^\lambda,\cH)$, if $\lambda(vg)f(g)=f(g)$ for almost all $g\in G$ and for all $v\in V_\lambda$, then $f=0$.

Now for any fixed $g\in G$, the condition $\lambda(vg)f(g)=f(g)$ for all $v\in V_\lambda$ implies that either $f(g)=0$ or $\lambda(vg)=1$ for all $v\in V_\lambda$. Writing $\lambda(vg)=e^{ivg\cdot  \alpha}$ we see that $\lambda(vg)=1$ for all $v\in V_\lambda$ if and only if $\alpha g^*\in V_\lambda^\bot=\R\alpha$. Next, noting that $Q^*(\alpha g^*)=Q^*(\alpha)$ for any $g\in G=\SO^+_{Q}(\R)$, if $\alpha g^*=c\alpha$ then $c^2Q^*(\alpha)=Q^*(c\alpha)=Q^*(\alpha)\neq 0$ so $c^2=1$, implying that $\alpha g^*=\pm \alpha$. Hence, if $f\in L^2(\tilde{G}/\tilde{G}^\lambda,\cH)$, satisfies $\lambda(vg)f(g)=f(g)$ for all $v\in V_\lambda$ then up to a null set, $f$ is supported on the set $\{g\in G/G^\lambda\ |\ \alpha g^*=\pm \alpha\}$ containing at most two points in $G/G^\lambda$. Since $f\in L^2(G/G^\lambda,\cH)$ is only defined up to its values on null sets, the only such element is the zero vector.
\end{proof}

The final ingredient for the induction argument is the following  argument going back to Burger and Sarnak \cite{BurgerSarnak1991}.
Let $G=\SO^+(n-1,1)$ and recall that for any unitary representation $\pi$ of $G$ not weakly containing the trivial representation, the parameter $\alpha(\pi)=\frac{n-2}{p(\pi)}$ characterizes the fastest decay rate of matrix coefficients of $\pi$ restricted to a fixed Cartan subgroup $A\leq G$ . To carry over this reduction argument it is more convenient to work with this parameter $\alpha(\pi)$.
Now inside $G$ we have a sequence of closed subgroups
$$G=G^{(1)}\supset G^{(2)}\supset G^{(n-2)}\supset A,$$ 
with $G^{(j)}\cong \SO^+(n-j,1)$ for $2\leq j\leq n-2$, all containing the same fixed Cartan group $A$ (so that for each $G^{(j)}$ we have a decomposition $G^{(j)}=K_jA^+K_j$ with $K_j\leq G^{(j)}$ a maximal compact subgroup). Since any $K$-finite vector in $\pi$ is also a $K_j$-finite vector in $\pi|_{G^{(j)}}$ and the parameter $\alpha(\pi)$ depends only on the $A$-action on $K$-finite vectors, we have that 
$\alpha(\pi)\geq \alpha(\pi|_{G^{(j)}})$. This leads to the following simple lemma reducing the proof of \thmref{t:rank1b} to studying the restriction representation $\pi|_{G^{(j)}}$ for some $1\leq j\leq n-2$.

\begin{Lem}\label{l:reduction}
Keep the notation and assumptions as above and let $\kappa\in [1,2]$.  If $p(\pi|_{G^{(j)}})\leq \kappa(n-j-1)$ for some $1\leq j\leq n-2$ then $p(\pi)\leq \kappa(n-2)$ (and $\alpha(\pi)\geq 1/\kappa$).
\end{Lem} 
\begin{proof}
Suppose $p(\pi|_{G^{(j)}})\leq \kappa(n-j-1)$ for some $2\leq j\leq n-2$, then by the relation \eqref{equ:reductionrelation} we have $\alpha(\pi)\geq\alpha(\pi|_{G^{(j)}})=\frac{n-j-1}{p(\pi|_{G^{(j)}})}\geq \frac1\kappa$. Again by \eqref{equ:reductionrelation} we have $p(\pi)=\frac{n-2}{\alpha(\pi)}\leq \kappa(n-2)$, finishing the proof.
\end{proof}

We can now give the 
\begin{proof}[Proof of \thmref{t:rank1b}]
The proof is by induction on $n$. When $n=4$, by our hypothesis we have $p(\tilde{G})_0= 2$, so we may assume $n\geq 5$.
%For $n=4$ by our hypothesis  %The basis of the induction, $n=4$, by our hypothesis $p(\tilde{G}) follows from \cite[Theorem 7.3.9]{Zimmer1984} so we may assume $n\geq 4$. 
Let $\tilde{\pi}$ denote a representation of $\tilde{G}=\SO^+_Q(\R)\ltimes \R^n$ with no nontrivial $\R^n$-invariant vectors and let $\pi=\tilde{\pi}|_G$.  Since almost every irreducible component of $\tilde\pi$ has no nontrivial $\R^n$-invariant vectors we may assume that $\tilde\pi$ is irreducible. 
Then $\tilde{\pi}=\ind_{\tilde G^\lambda}^{\tilde G}\sigma$ for some nontrivial unitary character $\lambda$ of $\R^n$, and an irreducible unitary representation $\sigma$ of $\tilde{G}^\lambda$.

Now from the discussion above,  either $G^\lambda$ is amenable, in which case $\pi$ is tempered, or $G^\lambda=\SO^+_{Q'}$ with $Q'$ of signature $(n-2,1)$. In the second case, let $\tilde{G}_2=G^\lambda\ltimes V_\lambda \cong \SO_{Q'}^+\ltimes \R^{n-1}$ with $V_\lambda=\ker(\lambda)$ as above. By \lemref{l:Vlambda} the restriction of $\tilde{\pi}$ to $\tilde{G}_2$ has no nontrivial $V_\lambda$-invariant vectors, 
so decomposing it as a direct integral 
$$\tilde{\pi}|_{\tilde{G}_2}=\int^\oplus \tilde\pi_x,$$
with $\tilde\pi_x$ irreducible, since $\tilde{\pi}|_{\tilde{G}_2}$ has no nontrivial $V_{\lambda}$-invariant vectors then $\tilde\pi_x$ has no nontrivial $V_\lambda$-invariant vectors for almost every $x$, and hence by induction $p(\tilde{\pi}_x|_{G^{\lambda}})\leq (n-3)$ for almost every $x$, implying that
$p(\tilde{\pi} |_{G^{\lambda}})\leq (n-3)$. Finally since $\tilde{\pi}|_{G^{\lambda}}=\pi|_{G^{\lambda}}$ by \lemref{l:reduction} we get $p(\pi)\leq (n-2)$ finishing the proof.
\end{proof}
\begin{rem}
Unconditionally, we can use the same argument starting with $n=3$ as the base of the induction (i.e. the aforementioned result of Kazhdan \cite{Kazhdan1967}), giving the bound $p(\tilde{G})_0\leq 2(n-2)$ which recovers the result of Wang \cite{Wang2014}.
\end{rem}

%\bibliographystyle{alpha}
%\bibliography{DKbibliog}

\end{document}